\documentclass[11pt]{elsarticle}

\usepackage{amssymb}
\usepackage{amsmath, amsthm}
\usepackage{mwe}

\usepackage{mathrsfs}
\usepackage{algorithm}
\usepackage{algpseudocode} 

\usepackage[
top=2.5cm,
bottom=2.5cm,
left=3cm,
right=3cm
]{geometry}

\usepackage{hyperref}
\usepackage{dsfont}
\usepackage{xcolor}
\usepackage{graphicx}
\usepackage{subcaption}

\newtheorem{theorem}{Theorem}[section]
\newtheorem{lemma}[theorem]{Lemma}
\newtheorem{corollary}[theorem]{Corollary}
\newtheorem{proposition}[theorem]{Proposition}
\newtheorem{definition}[theorem]{Definition}
\newtheorem{remark}{Remark}

\newtheorem{example}{Example}
\newtheorem{Notation}{Notation}
\numberwithin{equation}{section}

\journal{*****************}

\begin{document}

\begin{frontmatter}

\title{Controllability of finite-dimensional linear fractional systems under uncertain parameters} %% Article title

\author{Idriss Boutaayamou$^a$}
\ead{dsboutaayamou@gmail.com}
\author{Fouad Et-tahri$^{*,b}$}
\ead{f.ettahri@uca.ac.ma}
\author{Lahcen Maniar$^{b,c}$}
\ead{maniar@uca.ac.ma}

% Author affiliation
\affiliation{organization={Department of Mathematics},%Department and Organization
	addressline={Lab-SIV, Polydisciplinary Faculty-Ouarzazate, Ibnou Zohr University}, 
	city={Ouarzazate},
	postcode={45000}, 
	state={B.P. 638},
	country={Morocco}}
\affiliation{organization={Cadi Ayyad University, UCA, Faculty of Sciences Semlalia, Laboratory of Mathematics, Modeling and Automatic Systems},%Department and Organization
            addressline={B.P. 2390}, 
            city={Marrakesh},
            postcode={40000},
            country={Morocco}}
%\affiliation{organization={$^*$ Corresponding author. Faculty of Sciences-Agadir, %Lab-SIV, Ibnou Zohr University},%Department and Organization
 %           addressline={B.P. 8106}, 
  %          city={Agadir},
   %         postcode={80000},
    %        country={Morocco}}
\affiliation{organization={The UM6P Vanguard Center, Mohammed VI Polytechnic University},%Department and Organization
            addressline={Hay Moulay Rachid}, 
            city={Ben Guerir},
            postcode={43150},
            country={Morocco}}

%% Abstract
\begin{abstract}
This paper investigates the controllability of finite-dimensional linear fractional systems involving an uncertain parameter. We establish new results on the simultaneous and average controllability. In particular, we show that average controllability can be characterized by the so-called average Kalman rank condition and the average Gramian matrix. Moreover, using the average Gramian matrix, we design an open-loop control with minimal energy. These results can be seen as a natural generalization of the classical results known for systems with integer-order derivatives. Finally, numerical simulations are provided to robustly validate the theoretical findings, with a focus on the Rössler system for a damping parameter following either a Gaussian or a uniform distribution.

\end{abstract}

%\item **************
%
%\item **************
%
%\item ****************.

%% Keywords

\end{frontmatter}



\section{Introduction}
\par 
The controllability of dynamical systems using open-loop control is a central issue in many areas of science and engineering. In practical applications, the equations modeling such systems often involve uncertain parameters and exhibit memory effects. This has led to a growing interest in the study of control systems involving time-fractional derivatives and uncertain parameters, a field that has become particularly active in recent years.
Consider a parameter-dependent finite-dimensional linear fractional control system:
\begin{equation}\label{sys}
	\begin{cases}
		\partial^{\alpha}_t x_{\sigma}(t) =A(\sigma) x_{\sigma}(t)+B(\sigma)u(t),& t>0,\\
		x_{\sigma}(0)=x_{0,\sigma},
	\end{cases}
\end{equation}
where $\sigma\in\mathfrak{S}$ is a random parameter (representing the system's parameter) following a probability measure  $\mu$, with $(\mathfrak{S},\mathcal{F},\mu)$ a probability space, $H=\mathbb{R}^{n}$ is the state space and $U=\mathbb{R}^{m}$ the control one, $x_{\sigma}: [0,\infty)\rightarrow H$ is the state of the system, $\partial^{\alpha}_t$ denotes the left Caputo fractional derivative of order $\alpha\in (0,1)$ defined by
\begin{equation*}\label{cap01}
	\partial_t^\alpha x_\sigma(t):=\frac{1}{\Gamma(1-\alpha)} \int_0^t(t-s)^{-\alpha} x_\sigma^{\prime}(s)\,d s,
\end{equation*}
whenever the right-hand side is well-defined, and $\Gamma$ is the standard Gamma function.
\par 
For $\mu$-a.e, $\sigma\in\mathfrak{S}$, we assume that $A(\sigma)\in \mathbb{R}^{n\times n}$, $B(\sigma)\in\mathbb{R}^{n\times m}$,  $x_{0,\sigma}\in H$ and the input (control) $u=u(t)$ is independent of the random parameter $\sigma$. Given $x_{0,\sigma}\in H$ and $u\in L^{\infty}_{\text{loc}}(0,\infty; U)$, the Duhamel formula gives:
\begin{eqnarray*}
	x_{\sigma}(t;x_{0,\sigma}; u)=E_{\alpha}(t^{\alpha}A(\sigma))x_{0,\sigma} +\int_{0}^{t}(t-s)^{\alpha-1}E_{\alpha,\alpha}((t-s)^{\alpha}A(\sigma))B(\sigma)u(s)ds,\quad t\geq 0,
\end{eqnarray*}
where $E_\alpha$ and $E_{\alpha,\alpha}$ are the Mittag-Leffler matrix functions defined in \eqref{Mittag-Leffler one parameter} and \eqref{Mittag-Leffler two parameters}, respectively.
\par 
In this work, inspired by \cite{loheac2016averaged, zuazua2014averaged}, where the average and simultaneous controllability of finite-dimensional systems was studied, we aim to extend their results to the setting of fractional systems. 
\par 
Given a fixed time horizon $T>0$, the controllability of the system \eqref{sys} for almost every 
$\sigma\in\mathfrak{S}$ at time $T$ is characterized by the Kalman rank condition (see \cite[Corollary 2]{jolic2023}):
\begin{align*}
	\mbox{rank}\left[B(\sigma)\;A(\sigma)B(\sigma)\;\cdots\;A^{n-1}(\sigma)B(\sigma)\right]=n,\quad a.e, \sigma\in\mathfrak{S}.
\end{align*}
However, in this case, the control would depend on the random parameter, which is not desirable in practice. It is natural to control at least the average of system \eqref{sys}:
\begin{eqnarray*}
	\forall x_{0,\sigma}\in L^1(\mathfrak{S},H;\mu), \forall x_1\in H, \exists u\in L^\infty(0,T;U), \; \mathbb{E}(x_{\sigma}(T;x_{0,\sigma};u))=x_1
\end{eqnarray*}
which is characterized by the following averaged Kalman rank condition as shown in Theorem \ref{Average Kalman Rank}:
\begin{align*}
	\mbox{rank}\left[ \mathbb{E}(A^{k}(\sigma)B(\sigma))\colon k=0,1,\cdots \right]=n.
\end{align*}
In the special case where \(A(\sigma)\) and \(B(\sigma)\) are scalar multiples of fixed matrices \(A\) and \(B\), respectively, the averaged Kalman rank condition reduces to the classical Kalman rank condition, as stated in Theorem \ref{Average Kalman Rank simple}.\\
A key condition for average controllability is the invertibility of the averaged Gramian matrix, which enables the explicit construction of an open-loop control with minimal energy as established in Theorem \ref{open-loop control} and Proposition \ref{HUM Control a}.

In general, the difficulty of average controllability lies in the fact that the dynamics of the avergae cannot be directly derived from system \eqref{sys}, except when the matrix $A$ is independent of the random parameter.
\par 
In the context of parameter-dependent control systems, the ideal situation is simultaneous controllability. In this case, it means that:
\begin{eqnarray*}
	\forall x_{0,\sigma}, x_{1,\sigma}\in L^2(\mathfrak{S},H;\mu), \exists u\in L^\infty(0,T;U), \; x_{\sigma}(T;x_{0,\sigma};u)=x_{1,\sigma},\quad a.e, \sigma\in\mathfrak{S}.
\end{eqnarray*}
This means that one can find a control, independent of the random parameter, that successfully controls each system for every realization, which is highly challenging task. In the infinite-dimensional setting $L^2(\mathfrak{S},H;\mu)$, the system \eqref{sys} fails to be simultaneously controllable as shown in Proposition \ref{infinite dimensional}. In contrast, for finite parameter sets $\mathfrak{S}$, simultaneous controllability aligns with classical controllability via an associated extended system as established in Proposition \ref{finitely countable}. In the general case, no matrix criterion, such as the Kalman rank condition, is available, even for classical derivatives. However, this notion can still be studied using a variational framework, as in \cite{loheac2016averaged}.

For the sake of simplicity, we considered time-invariant matrices. However, the results obtained in this paper can be extended to linear time-varying fractional systems by introducing the state-transition matrix.

\subsection{Related literature}
\par 
The control of averaged properties for finite-dimensional systems was first addressed in the seminal work \cite{zuazua2014averaged}, where the author characterized average controllability using the Kalman rank condition. A pioneering study on both average and simultaneous controllability in the context of ordinary differential equations was presented in \cite{loheac2016averaged}. In that work, the authors investigated these two notions and investigated simultaneous controllability by deriving it from averaged controllability using a variational approach combined with a penalization method.

\par 
In the context of partial differential equations (PDEs), a pioneering study on average controllability was presented in \cite{lu2016averaged}. In that work, the authors investigated the problem of controlling the average behavior of solutions to PDEs, focusing on specific cases of the transport, heat, and Schrödinger equations where the diffusivity is modeled as a random variable. In \cite{coulson2019average} and \cite{barcena2021averaged}, the authors studied the controllability properties of the heat equation with random diffusivity. In \cite{barcena2026averaged}, the authors established the null controllability of the Schrödinger equation with a randomly varying diffusivity following a general absolutely continuous distribution. They also showed the lack of simultaneous and exact controllability for any absolutely continuous random variable.

\par 
In the context of controllability for finite-dimensional time-fractional systems, we refer to \cite{jolic2023, lamrani2025controllability}, where exact controllability is characterized using the Gramian matrix and the Kalman rank condition. In the first work, the authors consider systems involving Caputo derivatives with time-dependent matrices, while the second focuses on systems with tempered fractional derivatives and time-independent matrices.

To the best of our knowledge, there are currently no available results on average and simultaneous controllability for finite-dimensional fractional systems involving the Caputo derivative. This motivated us to establish new and meaningful results in this direction, based on a generalization of classical results from integer-order linear control theory.

\subsection{Contents} The manuscript is structured as follows. In Section \ref{Sec2}, we review some preliminary concepts from fractional calculus, probability theory, and linear algebra. Section \ref{Sec3} is devoted to the well-posedness analysis of the primal system \eqref{sys} and its adjoint, as well as a duality relationship between them. Section \ref{Sec4} focuses on the analysis of simultaneous and average controllability using a variational approach. The main results concern the characterization of average controllability through matrix criteria, specifically the Gramian matrix and the Kalman rank condition in the average sense. Section \ref{Sec5} provides applications illustrating the obtained results, with a focus on the Rössler system. Finally, Section \ref{Sec6} concludes the paper.
\section{Preliminaries} \label{Sec2}
In this section, we provide an overview of fractional calculus, probability, and linear algebra that will be needed throughout the remainder of the work.
\subsection{Overview of fractional calculus} \label{FDI}
In this section, we fix $T > 0$ and set $H = \mathbb{R}^n$ for some $n \geq 1$, equipped with its standard Euclidean structure, where $\langle\cdot,\cdot\rangle_{H}$ denotes its scalar product. The spaces $L^p(0,T;H)$ for $1 \le p \le \infty$ and $C([0,T];H)$ denote the usual functional spaces, while $W^{1,1}(0,T;H)$ stands for the space of absolutely continuous functions on $[0,T]$ with values in $H$.
Consider the following scalar function:
$$\omega_{\alpha}(t)=\begin{cases}
	\frac{1}{\Gamma(\alpha)t^{1-\alpha}} \quad &t>0\\
	0 & t\leq 0,
\end{cases}$$
where $\alpha>0$ and $\Gamma$ is the Gamma function.\\
The left and right Riemann–Liouville fractional integrals of order $\alpha>0$ are respectively defined as convolution product by the scalar function 
$\omega_\alpha\in L^1(0,T;\mathbb{R})$, as follows:
\begin{align*}
	(I^{\alpha}_{t}g)(t)&:=(\omega_{\alpha}\ast g)(t) =\int_{0}^{t}\omega_{\alpha}(t-s)g(s)ds, \\
	(I^{\alpha}_{T-t}g)(t)&:=(\omega_{\alpha}\ast g(T-\cdot))(T-t) =\int_{t}^{T}\omega_{\alpha}(s-t)g(s)ds,
\end{align*}
whenever the right-hand side is well-defined. Using the properties of the convolution product, it follows that \cite{jin2021fractional}:
\begin{lemma}[Continuity of fractional integrals] \label{Continuity of fractional integrals}
	Let $\alpha>0$ and $1\leq p\leq \infty$, then $I^{\alpha}_{t}$ (the same applies to $I^{\alpha}_{T-t}$) is continuous from $L^{p}(0,T;H)$ into itself
	and we have that:
	\begin{align*}
		\|I^{\alpha}_{t}g\|_{L^{p}(0,T;H)}\leq \frac{T^{\alpha}}{\Gamma(\alpha+1)}\|g\|_{L^{p}(0,T;H)}.
	\end{align*}
	It is also continuous from $W^{1,1}(0,T;H)$ into itself.
\end{lemma}
There are several types of fractional integration by parts formulas. For fractional
integrals, we have that:
\begin{lemma}[Fractional integration by parts \cite{jin2021fractional}] \label{Fractional integration by parts}
	For any $f\in L^p(0,T;H), g\in L^q(0,T;H)$, $p,q\geq 1$ and $\alpha>0$ with $p^{-1}+q^{-1}\leq 1+\alpha$, then $t\mapsto \langle (I_{t}^{\alpha}f)(t), g(t)\rangle_H$ and $t\mapsto \langle f(t), (I_{T-t}^{\alpha}g)(t)\rangle_H$ are in $L^{1}(0,T;\mathbb{R})$, and
	the following identity holds
	\begin{align*}
		\int_{0}^{T} \langle (I_{t}^{\alpha}f)(t), g(t)\rangle_H dt=\int_{0}^{T}\langle f(t), (I_{T-t}^{\alpha}g)(t)\rangle_H dt.
	\end{align*}
\end{lemma}

We now introduce the notions of Caputo and Riemann–Liouville fractional derivatives.
\begin{itemize}
	\item  We define the left and right Caputo fractional derivatives of order $\alpha\in (0,1)$, respectively, by
	\begin{align*}
		\partial^{\alpha}_{t}g(t)&
		=(I^{1-\alpha}_{t}g^{\prime})(t) \quad\mbox{and}\quad \partial^{\alpha}_{T-t}g(t)=(I^{1-\alpha}_{T-t}g^{\prime})(t)
	\end{align*}
	whenever these quantities exist. Note that, if $g\in W^{1,1}(0,T;H)$, then $\partial^{\alpha}_{t}g, \partial^{\alpha}_{T-t}g$ exist and $\partial^{\alpha}_{t}g, \partial^{\alpha}_{T-t}g\in L^1(0,T;H)$.
	\item We also define the left and right Riemann–Liouville fractional derivatives of order $\alpha\in (0,1)$, respectively, by
	\begin{align*}
		D^{\alpha}_{t}g(t)&=
		\frac{d}{dt}(I^{1-\alpha}_{t}g)(t) \quad\mbox{and}\quad D^{\alpha}_{T-t}g(t)=-\frac{d}{dt}(I^{1-\alpha}_{T-t}g)(t)
	\end{align*}
	whenever these quantities exist. Note that, if $g\in W^{1,1}(0,T;H)$, then $D^{\alpha}_{t}g, D^{\alpha}_{T-t}g$ exist and $D^{\alpha}_{t}g, D^{\alpha}_{T-t}g\in L^1(0,T;H)$.
\end{itemize}

The following integration formula is fundamental and allows the derivation of the adjoint system corresponding to the primal system involving the Riemann–Liouville fractional derivative.
\begin{proposition}
	For $\alpha\in (0,1)$ and $T>0$, assume that:
	\begin{itemize}
		\item $f\in C([0,T];H)$ such that $\partial^{\alpha}_tf \in L^{\infty}(0,T;H)$ and $f^{\prime}\in L^{q}(0,T;H)$ for somme  $q\geq\frac{1}{1-\alpha}$,
		\item $g\in L^{1}(0,T;H)$ such that $D^{\alpha}_{T-t}g\in L^1(0,T;H)$.
	\end{itemize}
	Then 
	\begin{align}
		\int_{0}^{T} \langle\partial^{\alpha}_{t}f(t), g(t)\rangle_H dt=\left[\langle f(t), (I^{1-\alpha}_{T-t}g)(t) \rangle_H\right]_{t=0}^{t=T} +\int_{0}^{T} \langle f(t), D^{\alpha}_{T-t}g(t)\rangle_H dt. \label{IBPF}
	\end{align}
\end{proposition}
\begin{proof}
	Based on the Caputo derivative definition, we rewrite
	\begin{align*}
		\int_{0}^{T} \langle\partial^{\alpha}_{t}f(t), g(t)\rangle_H dt=\int_{0}^{T} \left\langle (I^{1-\alpha}_{t}f^{\prime})(t), g(t)\right\rangle_H dt.
	\end{align*}
	Applying Lemma \ref{Fractional integration by parts} for $1-\alpha$ (instead of $\alpha$), $p=1$ and $q\geq \frac{1}{1-\alpha}$, we obtain
	\begin{align}
		\int_{0}^{T} \langle\partial^{\alpha}_{t}f(t), g(t)\rangle_H dt=\int_{0}^{T} \left\langle f^{\prime}(t), (I^{1-\alpha}_{T-t}g)(t) \right\rangle_H dt. \label{Integral}
	\end{align}
	We apply the standard integration by parts formula in $W^{1,1}(0,T;H)$ to 
	$f$ and $I^{1-\alpha}_{T-t}g$, the integral on the right-hand side of \eqref{Integral} leads to the conclusion \eqref{IBPF}.
\end{proof}
\begin{remark} \label{Rem}
	$(I^{1-\alpha}_{T-t}g)(t)$ is well defined as an element of $H$ for all $t\in [0,T]$, due to $I^{1-\alpha}_{T-t} g \in W^{1,1}(0,T;H)$, since $g \in L^1(0,T;H)$ (then $I^{1-\alpha}_{T-t} g\in L^1(0,T;H)$ see Lemma \ref{Continuity of fractional integrals}) and $\frac{d}{dt}(I^{1-\alpha}_{T-t}g)(t)=-D^{\alpha}_{T-t}g \in L^1(0,T;H)$.
\end{remark}
\subsection{Overview of Probability Concepts}
\begin{itemize}
	\item For any probability space $(\mathfrak{S}, \mathcal{F}, \mu)$ and any Banach space $E$ equipped with its Borel sigma-algebra, 
	$\mathcal{R}(\mathfrak{S}; E)$ denotes the set of random variables from $\mathfrak{S}$ to $E$, i.e., the measurable functions from $\mathfrak{S}$ to $E$. For $p\in [1,\infty]$, the space 
	$L^{p}(\mathfrak{S}, E; \mu)$ is defined by
	\begin{align*}
		L^{p}(\mathfrak{S}, E; \mu)&=\{f\in \mathcal{R}(\mathfrak{S}; E)\;:\; \int_{\mathfrak{S}}\|f(\sigma)\|^{p}_{E}d\mu(\sigma)<\infty\}, \; (1\le p<\infty),\\
		L^{\infty}(\mathfrak{S}, E; \mu)&=\{f \in \mathcal{R}(\mathfrak{S}; E) \;:\; \operatorname*{ess\,sup}_{\sigma \in \mathfrak{S}} \|f(\sigma)\|_E < \infty\}.
	\end{align*}
	\item Let $f\in L^{1}(\mathfrak{S}, E; \mu)$, the average  (or mathematical expectation) of $f$, denoted by $\mathbb{E}(f)$, is defined as:
	\begin{eqnarray*}
		\mathbb{E}(f):= \int_{\mathfrak{S}}f(\sigma)d\mu(\sigma).
	\end{eqnarray*}
	If $f\in L^{1}(\mathfrak{S}, H; \mu)$ and $y\in H$, then $\langle f, y\rangle_H\in L^{1}(\mathfrak{S}, \mathbb{R}; \mu)$
	and we have that:
	\begin{align*}
		\mathbb{E}\left(\langle f ,y\rangle_H\right)=\langle \mathbb{E}(f),y\rangle_H.
	\end{align*}
	\item For $p=2$ and $H$ is a Hilbert space, $L^{2}(\mathfrak{S}, H; \mu)$ is a Hilbert space endowed with the scalar product:
	\begin{align*}
		\langle f, g\rangle_{L^{2}(\mathfrak{S}, E; \mu)}:=&\int_{\mathfrak{S}}\langle f(\sigma), g(\sigma)\rangle_{H}d\mu(\sigma)\\
		=&\mathbb{E}\left(\langle f ,g\rangle_H\right).
	\end{align*}
	We have the following usual inequalities:
	\begin{align}
		\left|\mathbb{E}(\langle f,g\rangle_H)\right|&\leq (\mathbb{E}(\|f\|_H^2))^{\frac{1}{2}}(\mathbb{E}(\|g\|_H^2))^{\frac{1}{2}} \nonumber\\
		&\leq \varepsilon \mathbb{E}(\|f\|_H^2) + \frac{1}{4\varepsilon}\mathbb{E}(\|g\|_H^2)\quad \forall \varepsilon>0. \label{Young}
	\end{align}
\end{itemize}

\subsection{Overview of Linear Algebra Concepts}
We consider $H=\mathbb{R}^n$ and $U=\mathbb{R}^m$ for $n,m\geq 1$, endowed with their standard Euclidean structures. The inner product and the norm on $H$ are denoted by $\langle \cdot,\cdot\rangle_H$ and $\|\cdot\|_H$, respectively (the same notation is used for $U$). Throughout this paper, $\mathbb{R}^{n\times n}$ and $\mathbb{R}^{n\times m}$  denote the spaces of real $n\times n$ and $n\times m$ matrices, respectively. For $M\in\mathbb{R}^{n\times m}$, we denote respectively by $M^{\top}, \mbox{rank}(M), \mbox{Ker}(M)$ and $\mbox{Range}(M)$ the transpose, the rank, the kernel, and the range of $M$, respectively. Recall the following orthogonality relation
\begin{align}
	\mbox{Ker}(M^{\top})=[\mbox{Range}(M)]^{\top}, \label{orthogonality relation}
\end{align}
where $[\mbox{Range}(M)]^{\top}$ denotes the orthogonal complement of the space $\mbox{Range}(M)$.\\
A matrix $M\in\mathbb{R}^{n\times n}$ is symmetric if $M^{\top}=M$, it is said to be positive semidefinite if $\langle Mx,x\rangle_H \geq 0$ for all $x\in H$, and it is said to be positive definite if $\langle Mx,x\rangle_H > 0$ for all $x\in H\setminus\{0\}$. Recall that any symmetric positive semidefinite matrix is diagonalizable (by the spectral theorem), with non-negative eigenvalues.
\section{Well-Posedness and duality relation} \label{Sec3}
The study of fractional derivatives is closely related to the Mittag-Leffler functions, first introduced in \cite{ML03}.
\begin{definition}
	Let $\alpha, \beta>0$. The one- and two-parameter Mittag-Leffler matrix functions of a matrix $A\in\mathbb{R}^{n\times n}$ are defined as
	\begin{align}
		E_{\alpha}(A)&=\sum_{n=0}^{\infty}\frac{A^{n}}{\Gamma(\alpha n+1)}\in \mathbb{R}^{n\times n}, \label{Mittag-Leffler one parameter}\\
		E_{\alpha,\beta}(A)&=\sum_{n=0}^{\infty}\frac{A^{n}}{\Gamma(\alpha n+\beta)}\in \mathbb{R}^{n\times n}, \label{Mittag-Leffler two parameters}
	\end{align}
	where $\Gamma$ is the Gamma function.
\end{definition}
\begin{remark}
	Stirling's formula ensures the absolute convergence of the above series in $\mathbb{R}^{n\times n}$. 
\end{remark}
\begin{remark}
	For $\alpha=\beta=1$, we recover the classical exponential matrix $\mathrm{e}^{A}$. A fundamental difference between the matrices $e^{A}$ and $E_{\alpha}(A)$ is that, unlike $e^{A}$ which is always invertible and $(e^{tA})_{t\ge 0}$ satisfies the semigroup law, $E_{\alpha}(A)$ is not necessarily invertible since its eigenvalues are of the form $\mbox{sp}((E_{\alpha}(A))=\{E_{\alpha}(\lambda)\colon \lambda\in\mbox{sp}(A)\}$, where $E_{\alpha}(z)$ is the Mittag-Leffler function, which generally has roots (see \cite{gorenflo2020mittag, popov2013distribution}).
\end{remark}
The following result follows directly from \cite[Theorem 6]{Bo18}.
\begin{proposition} \label{Well-Posedness}
	Let $\alpha\in (0,1)$, $T>0$,  $x_{0,\sigma}\in H, A(\sigma)\in\mathbb{R}^{n\times n}$ and $B(\sigma)\in\mathbb{R}^{n\times m}$ for almost every $\sigma\in\mathfrak{S}$ and $u\in L^{\infty}(0,T;U)$. Then, for almost every $\sigma\in\mathfrak{S}$, the system \eqref{sys} has a unique solution $x_\sigma\in C([0,T];H)$ such that  $\partial^{\alpha}_{t}x_\sigma\in L^{\infty}(0,T;H)$. Moreover, the solution is given by the Duhamel formula:
	\begin{align}
		x_{\sigma}(t; x_{0,\sigma}, u)=E_{\alpha}(t^{\alpha}A(\sigma))x_{0,\sigma} +\int_{0}^{t}(t-s)^{\alpha-1}E_{\alpha,\alpha}((t-s)^{\alpha}A(\sigma))B(\sigma)u(s)ds,\quad t\in [0,T]. \label{Duhamel formula}
	\end{align}
\end{proposition}
The following corollary concerns the regularity of the solutions to \eqref{sys} with respect to the random parameter.
\begin{corollary} \label{Regula Random}
	Assume that $A\in L^{\infty}(\mathfrak{S},\mathbb{R}^{n\times n};\mu)$ and $B\in L^{\infty}(\mathfrak{S},\mathbb{R}^{n\times m};\mu)$. Then, for all $x_{0,\sigma}\in L^p(\mathfrak{S},H;\mu)$ for $p\in [1,\infty]$ and $u\in L^{\infty}(0,T;U)$, the unique solution of \eqref{sys}, given by \eqref{Duhamel formula}, satisfies $x_\sigma\in C([0,T];L^{p}(\mathfrak{S},H;\mu))$, and there exists a positive constant $c$ such that:
	\begin{align}
		\|x_\sigma\|_{C([0,T];L^{p}(\mathfrak{S},H;\mu))}\leq c(\|x_{0,\sigma}\|_{L^{p}(\mathfrak{S},H;\mu)}+ \|u\|_{L^{\infty}(0,T;U)}). \label{Energy}
	\end{align}
\end{corollary}
\begin{proof}
	Obviously, we have
	\begin{align}
		x_{\sigma}(t; x_{0,\sigma}, u) = x_{\sigma}(t; x_{0,\sigma}, 0) + x_{\sigma}(t; 0, u), \quad t \in [0,T],\; \text{a.e.},\; \sigma \in \mathfrak{S}. \label{E1}
	\end{align}
	Let $a=\|A(\sigma)\|_{L^{\infty}(\mathfrak{S},\mathbb{R}^{n\times n};\mu)}$ and $b=\|A(\sigma)\|_{L^{\infty}(\mathfrak{S},\mathbb{R}^{n\times m};\mu)}$. Using the Duhamel formula \eqref{Duhamel formula}, we obtain the estimate
	\begin{align*}
		\|x_\sigma(t; x_{0,\sigma}, 0)\|_{L^{p}(\mathfrak{S},H;\mu)}\leq E_{\alpha}(T^{\alpha}a) \|x_{0,\sigma}\|_{L^{p}(\mathfrak{S},H;\mu)},\quad t\in [0,T]. 
	\end{align*}
	As a consequence, the map $t\mapsto x_\sigma(t; x_{0,\sigma}, 0)$ belongs to $C([0,T]; L^{p}(\mathfrak{S},H;\mu))$, and we further have
	\begin{align}
		\|x_\sigma(\cdot; x_{0,\sigma}, 0)\|_{C([0,T]; L^{p}(\mathfrak{S},H;\mu))}\leq E_{\alpha}(T^{\alpha}a) \|x_{0,\sigma}\|_{L^{p}(\mathfrak{S},H;\mu)}. \label{E2}
	\end{align}
	Moreover, we have the estimate
	\begin{align*}
		\|x_{\sigma}(t; 0, u)\|_{H}\le \frac{T^{\alpha}E_{\alpha,\alpha}(T^{\alpha}a)b}{\alpha}  \|u\|_{L^\infty(0,T;U)},\quad t\in [0,T].
	\end{align*}
	Since $\mu$ is a probability measure, it follows that
	\begin{align*}
		\|x_{\sigma}(t; 0, u)\|_{L^p(\mathfrak{S},H; \mu)}\le \frac{T^{\alpha}E_{\alpha,\alpha}(T^{\alpha}a)b}{\alpha}  \|u\|_{L^\infty(0,T;U)},\quad t\in [0,T].
	\end{align*}
	Hence the map $t\mapsto x_\sigma(t;0; u)$ belongs to $C([0,T]; L^{p}(\mathfrak{S},H;\mu))$. This continuity follows from the fact that the convolution of an element in $L^1$
	with an element in $L^\infty$
	is continuous. Moreover, we have
	\begin{align}
		\|x_\sigma(\cdot; 0, u)\|_{C([0,T]; L^{p}(\mathfrak{S},H;\mu))}\leq \frac{T^{\alpha}E_{\alpha,\alpha}(T^{\alpha}a)b}{\alpha}  \|u\|_{L^\infty(0,T;U)}. \label{E3}
	\end{align}
	Finally, combining \eqref{E1}, \eqref{E2}, and \eqref{E3}, we obtain the desired estimate \eqref{Energy}.
\end{proof}
\begin{remark}
	The boundedness of the matrices $A$ and $B$ with respect to the randomness is used only to ensure the existence of the expectation of $x_\sigma(T;x_{0,\sigma};u)$. Our analysis, however, extends to more general pairs of matrices—referred to as admissible pairs—for which the solution at time $T$ is integrable with respect to the probability measure. Throughout the remainder of the paper, we assume that $A$ and $B$ are bounded with respect to the randomness.
\end{remark}
Controllability, as usual, is dual to an observability inequality for the adjoint system. 
Using the integration by parts formula \eqref{IBPF}, we obtain the adjoint system corresponding to \eqref{sys} as follows:
\begin{equation}\label{adj_sys}
	\begin{cases}
		D^{\alpha}_{T-t} y_\sigma(t) =A(\sigma)^{\top} y_\sigma(t),\; 0<t<T,\\
		(I^{1-\alpha}_{T-t}y_{\sigma})(t)|_{t=T}=y_{T,\sigma},
	\end{cases}
\end{equation}
where $D^{\alpha}_{T-t}$ (resp. $I^{1-\alpha}_{T-t}$) denotes the right Riemann-Liouville time fractional derivative of order $\alpha$
(resp. the right Riemann-Liouville time fractional integral of order $1-\alpha$) as introduced in Subsection \ref{FDI}.
\begin{proposition}\label{Well-Posedness_adj}
	Let $\alpha\in (0,1)$, $T>0$ and $y_{T,\sigma}\in H$. Then, for almost every $\sigma\in\mathfrak{S}$ the system \eqref{adj_sys} has a unique solution $y_\sigma\in L^{1}(0,T;H)$ such that $D^{\alpha}_{T-t} y_\sigma\in  L^{1}(0,T; H)$. Moreover, the solution is given by the Duhamel formula:
	\begin{align}
		y_\sigma(t; y_{T,\sigma})=(T-t)^{\alpha-1}E_{\alpha,\alpha}((T-t)^{\alpha}A(\sigma)^{\top})y_{T,\sigma},\quad t\in [0,T). \label{Duhamel formula_adj}
	\end{align}
\end{proposition}
The proof of Proposition \ref{Well-Posedness_adj} can be found in Theorem 4.2 of \cite{idczak2011existence}.
\begin{remark}
	Note that the final condition in \eqref{adj_sys} makes sense due to $I^{1-\alpha}_{T-t} y_{\sigma} \in W^{1,1}(0,T;H)$, as explained in Remark \ref{Rem}.
\end{remark}
We now present the regularity properties of the solutions to \eqref{adj_sys} with respect to the random parameter.
\begin{corollary} \label{Regula Random adj}
	Assume that $A\in L^{\infty}(\mathfrak{S},\mathbb{R}^{n\times n};\mu)$. Then, for all $y_{T,\sigma}\in L^p(\mathfrak{S},H;\mu)$ for $p\in [1,\infty]$, the unique solution of \eqref{adj_sys}, given by \eqref{Duhamel formula_adj}, satisfies $y_\sigma\in C([0,T);L^{p}(\mathfrak{S},H;\mu))$.
\end{corollary}
\begin{proof}
	Let $a=\|A(\sigma)\|_{L^{\infty}(\mathfrak{S},\mathbb{R}^{n\times n};\mu)}$. Using the Duhamel formula \eqref{Duhamel formula_adj}, we obtain the estimate
	\begin{align*}
		\|y_\sigma(t; y_{T,\sigma})\|_{L^{p}(\mathfrak{S},H;\mu)}\leq (T-t)^{\alpha-1}E_{\alpha,\alpha}(T^{\alpha}a) \|y_{T,\sigma}\|_{L^{p}(\mathfrak{S},H;\mu)}<\infty \quad \forall t\in [0,T). 
	\end{align*}
	Hence $y_\sigma\in C([0,T);L^{p}(\mathfrak{S},H;\mu))$.
\end{proof}
\begin{remark}\label{Reg Integ adj}
	Using termwise integration along with the classical identity linking the Beta and Gamma functions (see formula (A.34) in \cite{gorenflo2020mittag}), we obtain
	\begin{align}
		I^{1-\alpha}_{T-t} y_{\sigma}(t; y_{T,\sigma})=E_{\alpha}((T-t)^{\alpha}A(\sigma)^{\top})y_{T,\sigma}\quad \forall t\in [0,T]. \label{PIF}
	\end{align}
	In particular,
	$I^{1-\alpha}_{T-t} y_{\sigma}\in C([0,T];L^{p}(\mathfrak{S},H;\mu))$ if $y_{T,\sigma}\in L^{p}(\mathfrak{S},H;\mu)$.
\end{remark}
\begin{Notation}
	For almost every $\sigma \in \mathfrak{S}$, $x_{\sigma}(t;x_{0,\sigma}; u)$ (resp. $y_{\sigma}(t;y_{T,\sigma})$) denotes the unique solution of \eqref{sys} (resp. \eqref{adj_sys}), corresponding to the initial condition $x_{0,\sigma} \in H$ and the control $u \in L^{\infty}(0,T;U)$ (resp. the final condition $y_{T,\sigma} \in H$).
\end{Notation}
The following duality relation between systems \eqref{sys} and \eqref{adj_sys} will be applied in the context of simultaneous controllability for $(p,q)=2$, and of average controllability for $(p,q)=(1,\infty)$.
\begin{lemma}
	Let $p,q\in [1,\infty]$ are conjugate, i.e., $\frac{1}{p}+\frac{1}{q}=1$, then for all $x_{0,\sigma}\in L^p(\mathfrak{S},H;\mu), y_{T,\sigma}\in L^q(\mathfrak{S},H;\mu)$ and $u\in L^{\infty}(0,T;U)$, we have that:
	\begin{align}
		\int_{0}^{T}\langle u(t),\mathbb{E}(B(\sigma)^{\top}y_\sigma(t;y_{T,\sigma})) \rangle_U dt=\mathbb{E}(\langle x_\sigma(T;x_{0,\sigma};u), y_{T,\sigma} \rangle_H) - \mathbb{E}(\langle x_{0,\sigma} , I^{1-\alpha}_{T-t}y_\sigma(0;y_{T,\sigma})\rangle_H). \label{relation duality}
	\end{align}
\end{lemma}
\begin{proof}
	For almost every 
	$\sigma\in\mathfrak{S}$, $x_{\sigma}(t;x_{0,\sigma};u)$ and $y_{\sigma}(t;y_{T,\sigma})$ are the solutions to \eqref{sys} and \eqref{adj_sys}, respectively, it follows that:
	\begin{align*}
		&\int_{0}^{T}\langle u(t),B(\sigma)^{\top}y_{\sigma}(t;y_{T,\sigma}) \rangle_U dt\\
		& =\int_{0}^{T}\langle \partial^{\alpha}_{t}x_{\sigma}(t; x_{0,\sigma};u),y_{\sigma}(t;y_{T,\sigma}) \rangle_H dt-\int_{0}^{T}\langle A(\sigma)x_{\sigma}(t; x_{0,\sigma};u),y_{\sigma}(t;y_{T,\sigma}) \rangle_H dt\\
		& =\int_{0}^{T}\langle \partial^{\alpha}_{t}x_{\sigma}(t; x_{0,\sigma};u),y_{\sigma}(t;y_{T,\sigma}) \rangle_H dt-\int_{0}^{T}\langle x_{\sigma}(t; x_{0,\sigma};u), A(\sigma)^{\top}y_{\sigma}(t;y_{T,\sigma}) \rangle_H dt\\
		& =\int_{0}^{T}\langle \partial^{\alpha}_{t}x_{\sigma}(t; x_{0,\sigma};u),y_{\sigma}(t; y_{T,\sigma}) \rangle_H dt-\int_{0}^{T}\langle x_{\sigma}(t; x_{0,\sigma};u), D^{\alpha}_{T-t}y_{\sigma}(t;y_{T,\sigma}) \rangle_H dt.
	\end{align*}
	Using the time-regularity of the solutions given in Propositions \ref{Well-Posedness} and \ref{Well-Posedness_adj}, we can apply the fractional integration by parts formula \eqref{IBPF}, then
	\begin{align}
		\int_{0}^{T}\langle u(t),B(\sigma)^{\top}y_{\sigma}(t;y_{T,\sigma}) \rangle_U dt
		& =\left[\langle x_{\sigma}(t; x_{0,\sigma};u),I^{1-\alpha}_{T-t}y_{\sigma}(t;y_{T,\sigma}) \rangle_H \right]_{t=0}^{t=T}. \label{lrd}
	\end{align}
	Finally, using the randomness-regularity of the solutions given in Corollaries \ref{Regula Random} and \ref{Regula Random adj}, and Remark \ref{Reg Integ adj}, we deduce the duality relation \eqref{relation duality} from \eqref{lrd} by taking expectations.
\end{proof}
\section{Controllability Analysis} \label{Sec4}
In this section, we prove several results regarding the controllability of system \eqref{sys}.
\subsection{Simultaneous controllability}
\begin{Notation}
	Let $T>0$ and define, for any initial state $x_{0,\sigma}\in L^2(\mathfrak{S},H;\mu)$, the set of the reachable states:
	\begin{align*}
		R(T; x_{0,\sigma}):=\{x_\sigma(T;x_{0,\sigma}; u)\colon u\in L^{\infty}(0,T;U)\},
	\end{align*}
	where, for almost every $\sigma \in \mathfrak{S}$, $x_\sigma(\cdot;x_{0,\sigma}; u)$ denotes the unique solution of \eqref{sys} associated with the initial data $x_{0,\sigma}$ and the control $u$.
\end{Notation}
\begin{remark}
	Proposition \ref{Well-Posedness} implies that $R(T; x_{0,\sigma})$ is well-defined and is an affine subspace of $L^2(\mathfrak{S},H;\mu)$, since the mapping $u\mapsto x_\sigma(T;x_{0,\sigma}; u)$ is affine from $L^{\infty}(0,T;U)$ to $L^2(\mathfrak{S},H;\mu)$.
\end{remark}
There are different notions of controllability in the simultaneous sense that need to be distinguished.
\begin{definition} \label{notion_simultaneous controllability}
	Let $T>0$.
	\begin{itemize}
		\item System \eqref{sys} is exactly simultaneous controllable in time T
		if, for every initial state $x_{0,\sigma}\in L^2(\mathfrak{S},H;\mu)$, the set of the reachable states
		$R(T, x_{0,\sigma})$ coincides with $L^2(\mathfrak{S},H;\mu)$. That is, for every $x_{0,\sigma}, x_{1,\sigma}\in L^2(\mathfrak{S},H;\mu)$, there exists $u\in L^{\infty}(0,T;U)$ such that:
		\begin{align}
			x_\sigma(T;x_{0,\sigma};u)=x_{1,\sigma},\; a.e,\;\sigma\in\mathfrak{S}. \label{exactly simultaneously controllable}
		\end{align}
		\item System \eqref{sys} is approximately simultaneous controllable in time T
		if, for every initial state $x_{0,\sigma}\in L^2(\mathfrak{S},H;\mu)$, the set of the reachable states
		$R(T, x_{0,\sigma})$ is dense in $L^2(\mathfrak{S},H;\mu)$. That is, for every $x_{0,\sigma}, x_{1,\sigma}\in L^2(\mathfrak{S},H;\mu)$ and every $\varepsilon>0$, there exists $u\in L^{\infty}(0,T;U)$ such that:
		\begin{align}
			\mathbb{E}(\|x_\sigma(T;x_{0,\sigma};u)-x_{1,\sigma}\|_{H}^2)<\varepsilon. \label{approximate simultaneously controllable}
		\end{align}
	\end{itemize}
\end{definition}
\begin{remark}
	It is clear that exact simultaneous controllability implies approximate simultaneous controllability.
\end{remark}
\begin{remark}
	In finite-dimensional spaces, affine subspaces are closed, so exact and approximate simultaneous controllability coincide when $L^2(\mathfrak{S},H;\mu)$ is finite-dimensional space. In the infinite-dimensional setting, the lack of simultaneous controllability was established in Proposition \ref{infinite dimensional}. Therefore, approximate simultaneous controllability would be a suitable alternative in this case. However, this notion is not the focus of the present paper and may be considered in future work.
\end{remark}
We begin with the simplest case where 
$\mathfrak{S}$ is finitely countable. In this case, $L^2(\mathfrak{S},H;\mu)$ is a finite-dimensional space, so exact and approximate simultaneous controllability coincide, and both are characterized by the classical notion of exact controllability (without uncertain parameters) for the associated extended system.
\begin{proposition} \label{finitely countable}
	Assume that $\mathfrak{S}=\{\sigma_1, \sigma_2, \cdots, \sigma_p\}$ is finitely countable. Then, the system \eqref{sys} is exactly simultaneously controllable in time $T$ if and only if
	\begin{align}
		\text{rank}[\mathbf{B}\;\;\mathbf{AB}\;\;\cdots\;\;\mathbf{A}^{p n-1}\mathbf{B}]=pn, \label{Kalman rank extended}
	\end{align}
	where $p=|\mathfrak{S}|$, $n=\dim(H)$, $\mathbf{A}=\text{diag}(A(\sigma_1),\cdots, A(\sigma_p))$ and $\mathbf{B}=\begin{pmatrix}B(\sigma_1)\\ \vdots \\ B(\sigma_p)\end{pmatrix}$.
\end{proposition}
\begin{proof}
	Obviously, the exact simultaneous controllability of system \eqref{sys} is equivalent to the classical notion of exact controllability for the extended system:
	\begin{equation}\label{extended system}
		\begin{cases}
			\partial^{\alpha}_t \mathbf{x}(t) =\mathbf{A} \mathbf{x}(t)+\mathbf{B}u(t),& 0<t<T,\\
			\mathbf{x}(0)=\mathbf{x}_{0},
		\end{cases}
	\end{equation}
	where $\mathbf{x}=\begin{pmatrix}
		x_{\sigma_1} \\ \vdots \\ x_{\sigma_p}
	\end{pmatrix}$, $\mathbf{x}_0=\begin{pmatrix} x_{0,\sigma_1}\\ \vdots \\ x_{0,\sigma_p}\end{pmatrix}$, $\mathbf{A}=\text{diag}(A(\sigma_1),\cdots, A(\sigma_p))$ and $\mathbf{B}=\begin{pmatrix}B(\sigma_1)\\ \vdots \\ B(\sigma_p)\end{pmatrix}$.
	The Kalman rank condition \eqref{Kalman rank extended} provides a characterization of the exact controllability for finite-dimensional time-fractional system \eqref{extended system} (see \cite[Corollary 2]{jolic2023}).
\end{proof}
In the case where $\mathfrak{S}$ is infinite, two scenarios are possible: the space $L^2(\mathfrak{S},H;\mu)$ may be finite- or infinite dimensional, depending on the probability measure $\mu$ (since $H$ is finite-dimensional). In the scenario where 
$L^2(\mathfrak{S},H;\mu)$ is infinite dimensional, and inspired by \cite[Theorem 1.2]{triggiani1977}, we prove the lack of exact simultaneous controllability. For that, we consider the weighted Hilbert space $L_{\alpha-1}^2(0,T,U)$, endowed with the scalar product:
\begin{align}
	\langle u, v\rangle_{L^2_{\alpha-1}(0,T;U)}&=\int_{0}^{T}(T-s)^{\alpha-1}\langle u(s),v(s) \rangle_U ds. \label{WHS}
\end{align}
This space is also used to characterize the unique control of minimal energy (see Proposition \ref{HUM Control s}). We now introduce the 
control-to-final-state operator $\mathfrak{L}_{T}: L_{\alpha-1}^2(0,T,U)\rightarrow L^{2}(\mathfrak{S},H;\mu)$, defined as 
\begin{align*}
	\mathfrak{L}_{T}(u)(\sigma)&:=\int_{0}^{T}(T-s)^{\alpha-1}E_{\alpha,\alpha}((T-s)^{\alpha}A(\sigma))B(\sigma)u(s)ds\\
	&=\int_{0}^{T}\mathcal{K}_{T}(\sigma,s)u(s)ds,
\end{align*}
where the Kernel $\mathcal{K}_T$ is given by:
\begin{align*}
	\mathcal{K}_T(\sigma,s)=(T-s)^{\alpha-1}E_{\alpha,\alpha}((T-s)^{\alpha}A(\sigma))B(\sigma)\in\mathbb{R}^{n\times m}.
\end{align*}
To show that $\mathfrak{L}_T$ is a compact operator, it suffices to prove that it is Hilbert–Schmidt, that is, 
\begin{align*}
	\sum_{m=0}^{\infty}\|\mathfrak{L}_{T}(u_m)\|^2_{L^2(\mathfrak{S},H;\mu)}<\infty,
\end{align*}
where $\{u_m\}$ a Hilbert basis of $L^2_{\alpha-1}(0,T; U)$.
\begin{lemma}
	$\mathfrak{L}_{T}: L^2_{\alpha-1}(0,T;U)\rightarrow L^{2}(\mathfrak{S},H;\mu)$ is a Hilbert-Schmidt operator.
\end{lemma}
\begin{proof}
	Let $\{u_m\}$ be a Hilbert basis of $L^2_{\alpha-1}(0,T; U)$ and denote by $\mathcal{K}^i_T(\sigma, s)$ the $i$-th row of the matrix $\mathcal{K}_T(\sigma, s)$.
	\begin{align*}
		\|\mathfrak{L}_{T}(u_m)(\sigma)\|_H^2
		&=\sum_{i=1}^n \left(\int_{0}^{T}\langle \mathcal{K}^i_{T}(\sigma,s), u_m(s)\rangle_U ds\right)^2\\
		&=\sum_{i=1}^n \langle \mathcal{K}^i_{T}(\sigma,s), u_m(s)\rangle^2_{L^2(0,T;U)}\\
		&=\sum_{i=1}^n \langle (T-s)^{1-\alpha}\mathcal{K}^i_{T}(\sigma,s), u_m(s)\rangle^2_{L^2_{\alpha-1}(0,T;U)}.
	\end{align*}
	Then 
	\begin{align*}
		\|\mathfrak{L}_{T}(u_m)\|^2_{L^2(\mathfrak{S},H;\mu)}&=\int_{\mathfrak{S}}\|\mathfrak{L}_{T}(u_m)(\sigma)\|_H^2 d\mu(\sigma)\\
		&=\int_{\mathfrak{S}}\sum_{i=1}^n \langle (T-s)^{1-\alpha}\mathcal{K}^i_{T}(\sigma,s), u_m(s)\rangle^2_{L^2_{\alpha-1}(0,T;U)}d\mu(\sigma).
	\end{align*}
	Therefore, applying Parseval's identity, we obtain
	\begin{align*}
		\sum_{m=0}^{\infty}\|\mathfrak{L}_{T}(u_m)\|^2_{L^2(\mathfrak{S},H;\mu)}
		&=\int_{\mathfrak{S}}\sum_{i=1}^n \sum_{m=0}^{\infty}\langle (T-s)^{1-\alpha}\mathcal{K}^i_{T}(\sigma,s), u_m(s)\rangle^2_{L^2_{\alpha-1}(0,T;U)}d\mu(\sigma)\\
		&=\int_{\mathfrak{S}}\sum_{i=1}^n \|(T-s)^{1-\alpha}\mathcal{K}^i_{T}(\sigma,s)\|^2_{L^2_{\alpha-1}(0,T;U)}d\mu(\sigma)\\
		&\leq \int_{\mathfrak{S}} \|(T-s)^{1-\alpha}\mathcal{K}_{T}(\sigma,s)\|^2_{L^2_{\alpha-1}(0,T;\mathcal{L}(H,U))}d\mu(\sigma)\\
		&=\int_{\mathfrak{S}}\int_{0}^{T} (T-s)^{(1-\alpha)}\|\mathcal{K}_{T}(\sigma,s)\|^2_{\mathcal{L}(H,U)}ds d\mu(\sigma)<\infty.
	\end{align*}
\end{proof}
We now state the result regarding the lack of simultaneous controllability.
\begin{proposition}\label{infinite dimensional}
	Assume that $L^2(\mathfrak{S},H;\mu)$ is infinite dimensional. Then, the system \eqref{sys} is not exactly simultaneous controllable in time $T$.
\end{proposition}
\begin{proof}
	Assume, by contradiction, that system \eqref{sys} is exactly simultaneous controllable in time $T$. Given that $L^\infty(0,T;U)\subset L_{\alpha-1}^2(0,T;U)$, then $\mathfrak{L}_T$ is surjective. Consequently
	\begin{align*}
		L^{2}(\mathfrak{S},H;\mu)=\bigcup_{n\geq 1}\overline{\mathfrak{L}_T(B_{n})}:=\bigcup_{n\geq 1} F_n,
	\end{align*}
	where $B_{n}$ denotes the closed ball of $L^2_{\alpha-1}(0,T;U)$ centered at $0$ with radius $n$.
	Since $\mathfrak{L}_T$ is compact, then $F_n$ is compact for any $n\geq 1$. Using the Baire Category Theorem (see \cite[Remark 1]{Br11}), there exists some $n_0\geq 1$ such that $\mbox{Int}(F_{n_0})\neq \varnothing$. Thus, $F_{n_0}$ contains a closed ball of $L^{2}(\mathfrak{S},H;\mu)$, which is compact (in the strong topology), it follows from Riesz's Theorem (see \cite[Theorem 6.5]{Br11}) that $L^{2}(\mathfrak{S},H;\mu)$ is finite-dimensional.
\end{proof}
Now, we study the general case by means of a variational approach. Let us first derive a necessary and sufficient condition for the exact simultaneous controllability of \eqref{sys}. The main idea is to characterize \eqref{exactly simultaneously controllable} through orthogonality in the space $L^2(\mathfrak{S},H;\mu)$:
\begin{align*}
	\mathbb{E}(\langle x_\sigma(T;x_{0,\sigma};u), y_{T,\sigma}\rangle_{H})=\mathbb{E}(\langle x_{1,\sigma}, y_{T,\sigma}\rangle_{H}) \quad \forall y_{T,\sigma}\in L^2(\mathfrak{S},H;\mu).
\end{align*}
Consequently, the duality relation \eqref{relation duality} yields the following characterization.
\begin{lemma} \label{Lemma: control}
	The system \eqref{sys} is driven from the initial condition $x_{0,\sigma} \in L^{2}(\mathfrak{S},H;\mu)$ to the final state $x_{1,\sigma} \in L^{2}(\mathfrak{S},H;\mu)$ in time $T$ in the sense of simultaneous controllability by a control $u \in L^{\infty}(0,T;U)$, i.e., $x_\sigma(T;x_0;u)=x_{1,\sigma},\;\text{a.e},\;\sigma\in\frak{S}$, if and only if
	\begin{align}
		\int_{0}^{T}\langle u(t),\mathbb{E}(B(\sigma)^{\top}y_\sigma(t;y_{T,\sigma})) \rangle_U dt=\mathbb{E}(\langle x_{1,\sigma}, y_{T,\sigma} \rangle_{H}) -\mathbb{E}(\langle x_{0,\sigma} , I^{1-\alpha}_{T-t}y_\sigma(0;y_{T,\sigma})\rangle_H), \label{First condition s}
	\end{align}
	for all $y_{T,\sigma}\in L^2(\mathfrak{S},H;\mu)$, 
	where $y_\sigma(t;y_{T,\sigma})$ denotes the unique solution of \eqref{adj_sys} associated with the final data $y_{T,\sigma}$.
\end{lemma}
Relation \eqref{First condition s} may be interpreted as an optimality condition of the critical points of the functional $\mathcal{J}^{s}: L^2(\mathfrak{S},H;\mu)\longrightarrow\mathbb{R}$:
\begin{align*}
	\mathcal{J}^{s}(y_{T,\sigma})=&\frac{1}{2}\int_{0}^{T}\|\mathbb{E}(B(\sigma)^{\top}y_{\sigma}(t;y_{T,\sigma}))\|^2_U \rho_{\alpha}(t)dt-\mathbb{E}(\langle x_{1,\sigma}, y_{T,\sigma} \rangle_{H}) \\
	& +\mathbb{E}(\langle x_{0,\sigma} , I^{1-\alpha}_{T-t}y_\sigma(0;y_{T,\sigma})\rangle_H),
\end{align*}
where $\rho_{\alpha}(t):=(T-t)^{1-\alpha},\; t\in [0,T]$.
\begin{remark}
	The presence of the weight $\rho_{\alpha}$ ensures integrability for all $\alpha \in (0,1)$. Indeed, based on the Duhamel formula \eqref{Duhamel formula_adj}, we obtain $\mathbb{E}(B(\sigma)^{\top}y_{\sigma}(t;y_{T,\sigma}))\underset{t\nearrow T}{\sim} c(T-t)^{\alpha-1}$, which implies $\mathbb{E}(B(\sigma)^{\top}y_{\sigma}(t;y_{T,\sigma}))\in L^{2}(0,T;U)$ if and only if $\alpha>\frac{1}{2}$.
\end{remark}
\begin{proposition} \label{Minimizer controllability}
	Let $x_{0,\sigma}, x_{1,\sigma}\in L^2(\mathfrak{S},H;\mu)$ and suppose that $\hat{y}_{T,\sigma}\in L^2(\mathfrak{S},H;\mu)$ is a minimizer of $\mathcal{J}^s$. Then, for the control 
	\begin{align}
		\hat{u}(t)
		&=\mathbb{E}\left(B(\sigma)^{\top}E_{\alpha,\alpha}((T-t)^{\alpha}A(\sigma)^{\top})\hat{y}_{T,\sigma}\right) \label{control s}
	\end{align}
	we have that:
	\begin{align}
		x_\sigma(T;x_{0,\sigma}; \hat{u})=x_{1,\sigma},\quad a.e,\; \sigma\in\mathfrak{S}. \label{ESC}
	\end{align}
\end{proposition}
\begin{proof}
	If $\hat{y}_{T,\sigma}$ is a minimizer of $\mathcal{J}^s$, then 
	\begin{align*}
		\lim_{h\to 0}\frac{\mathcal{J}^s(\hat{y}_{T,\sigma}+hy_{T,\sigma})-\mathcal{J}^s(\hat{y}_{T,\sigma})}{h}=0,
	\end{align*}
	for all $y_{T,\sigma}\in L^2(\mathfrak{S},H;\mu)$.This is equivalent to
	\begin{align*}
		&\int_{0}^{T}\langle \mathbb{E}(B(\sigma)^{\top}y_\sigma(t;\hat{y}_{T,\sigma})),\mathbb{E}(B(\sigma)^{\top}y_\sigma(t;y_{T,\sigma})) \rangle_U \rho_{\alpha}(t)dt-\mathbb{E}(\langle x_{1,\sigma}, y_{T,\sigma} \rangle_{H})\\ &+\mathbb{E}(\langle x_{0,\sigma} , I^{1-\alpha}_{T-t}y_\sigma(0;y_{T,\sigma})\rangle_H)=0,
	\end{align*}
	for all $y_{T,\sigma}\in L^2(\mathfrak{S},H;\mu)$, which, in view of Lemma \ref{Lemma: control}, implies that
	$$\hat{u}(t)
	=\rho_{\alpha}(t)\mathbb{E}(B(\sigma)^{\top}y_\sigma(t;\hat{y}_{T,\sigma}))\\
	=\mathbb{E}\left(E_{\alpha,\alpha}((T-t)^{\alpha}A(\sigma))\hat{y}_{T,\sigma}\right)\in L^\infty(0,T;U)$$
	is a control that satisfies \eqref{ESC}.
\end{proof}
\begin{theorem}
	The following statements are equivalent:
	\begin{itemize}
		\item System \eqref{sys} fulfills the exact simultaneous controllability.
		\item System \eqref{adj_sys} fulfills the
		exact simultaneous observability:
		\begin{align}
			\exists C_{obs}>0,\; \mathbb{E}(\|y_{T,\sigma}\|^{2}_H)\leq C_{obs} \int_{0}^{T}(T-t)^{1-\alpha}\| \mathbb{E}(B(\sigma)^{\top}y_{\sigma}(t;y_{T,\sigma}))\|^{2}_{U}dt, \label{exa_simu_observa_inequality}
		\end{align}
		for all $y_{T,\sigma}\in L^2(\mathfrak{S},H;\mu)$.
	\end{itemize}
\end{theorem}
\begin{proof}
	Let us suppose that \eqref{sys} is exactly simultaneous controllable and \eqref{adj_sys} is not exactly simultaneous observable. Then, for every
	$C>0$ there exists $y_{T,\sigma}\in L^2(\mathfrak{S},H;\mu)$ such that \eqref{exa_simu_observa_inequality} does not hold, and we can define a sequence  
	of functions $\{y_{T,\sigma}^k\}\subset L^2(\mathfrak{S},H;\mu)$, such that $\mathbb{E}(\|y_{T,\sigma}^k\|^2_{H})=1$ 
	and 
	\begin{align}
		\int_{0}^{T}(T-t)^{1-\alpha}\| \mathbb{E}(B(\sigma)^{\top}y_{\sigma}(t;y_{T,\sigma}^k))\|^{2}_{U}dt<\frac{1}{k}. \label{formulee}
	\end{align}
	Proposition \ref{infinite dimensional} ensures that $L^2(\mathfrak{S}, H; \mu)$ is finite-dimensional, and since the sequence $\{y_{T,\sigma}^k\}$ is bounded in $L^2(\mathfrak{S}, H; \mu)$, it admits a convergent subsequence (still denoted by the same name) that converges to $y^*_{T,\sigma}$ in $L^2(\mathfrak{S}, H; \mu)$.\\
	On the other hand, using Duhamel's formula \eqref{Duhamel formula_adj}, and noting that $A$ and $B$ are bounded, there exists a constant $C_T > 0$ such that
	\begin{align*}
		\| \mathbb{E}(B(\sigma)^{\top}y_{\sigma}(t;y_{T,\sigma}))\|^{2}_{U}&=(T-t)^{2(\alpha-1)}\| \mathbb{E}(B(\sigma)^{\top}E_{\alpha,\alpha}((T-t)^{\alpha}A(\sigma))y_{T,\sigma})\|^{2}_{U}\\
		&\leq (T-t)^{2(\alpha-1)}\mathbb{E}(\|B(\sigma)^{\top}E_{\alpha,\alpha}((T-t)^{\alpha}A(\sigma))y_{T,\sigma}\|^{2}_{U})\\
		&\leq C_T (T-t)^{2(\alpha-1)}\mathbb{E}(\|y_{T,\sigma}\|_H^2)
	\end{align*}
	for all $y_{T,\sigma}\in L^2(\mathfrak{S},H;\mu)$ and a.e, $t\in (0,T)$.
	Consequently, for a.e,\; $t\in (0,T)$:
	\begin{align*}
		(T-t)^{1-\alpha}\| \mathbb{E}(B(\sigma)^{\top}y_{\sigma}(t;y_{T,\sigma}^k))\|^{2}_{U}\rightarrow (T-t)^{1-\alpha}\| \mathbb{E}(B(\sigma)^{\top}y_{\sigma}(t;y_{T,\sigma}))\|^{2}_{U},\;\mbox{as}\; k\to \infty
	\end{align*}
	and 
	\begin{align*}
		(T-t)^{1-\alpha}\| \mathbb{E}(B(\sigma)^{\top}y_{\sigma}(t;y_{T,\sigma}^k))\|^{2}_{U}\leq C_T (T-t)^{\alpha-1}.
	\end{align*}
	Then, applying the dominated convergence theorem to \eqref{formulee}, we deduce that 
	\begin{align*}
		\int_{0}^{T}(T-t)^{1-\alpha}\| \mathbb{E}(B(\sigma)^{\top}y_{\sigma}(t;y^*_{T,\sigma}))\|^{2}_{U}dt=0.
	\end{align*}
	Then 
	\begin{align}
		\mathbb{E}(B(\sigma)^{\top}y_{\sigma}(t;y^*_{T,\sigma}))=0,\quad a.e, \; t\in (0,T).
	\end{align}
	Let $y_{T,\sigma}\in L^2(\mathfrak{S},H;\mu)$ be an arbitrary final state, setting $x_{0,\sigma}=0$ and $x_{1,\sigma}=y_{T,\sigma}$, there exists a control $u\in L^\infty(0,T;U)$ that satisfies $x_\sigma(T;0;u)=x_{1,\sigma}$, $\mu$-a.e, $\sigma\in\mathfrak{S}$. Applying \eqref{relation duality} to this control, we obtain 
	\begin{align*}
		\mathbb{E}(\|y_{T,\sigma} \|^{2}_H)=\int_{0}^{T}\langle u(t),\mathbb{E}(B(\sigma)^{\top}y_{\sigma}(t;y_{T,\sigma})) \rangle_U dt, 
	\end{align*}
	for all $y_{T,\sigma}\in L^2(\mathfrak{S},H;\mu)$. In particular,  
	\begin{align*}
		\mathbb{E}(\|y^*_{T,\sigma} \|^{2}_H)=\int_{0}^{T}\langle u(t),\mathbb{E}(B(\sigma)^{\top}y_{\sigma}(t;y^*_{T,\sigma})) \rangle_U dt=0.
	\end{align*}
	This is a contradiction with $\mathbb{E}(\|y^*_{T,\sigma} \|^{2}_H)=1$.\\ 
	\paragraph{}
	Conversely, let us suppose that the observability estimate \eqref{exa_simu_observa_inequality} is satisfied. Obviously, the functional $\mathcal{J}^s$ is convex and continuous in $L^2(\mathfrak{S},H;\mu)$ (for the fractional integral part, we use \eqref{PIF}). Thus, the
	existence of a minimum is ensured if $\mathcal{J}^s$ is coercive, i.e,
	\begin{align*}
		\mathcal{J}^s(y_{T,\sigma})\to \infty\;\mbox{when}\; \mathbb{E}(\|y_{T,\sigma}\|^2_H)\to \infty.
	\end{align*}
	The coercivity of the functional $\mathcal{J}^s$ follows immediately from the observability estimate \eqref{exa_simu_observa_inequality} and Remark \ref{Reg Integ adj}. Indeed,
	\begin{align*}
		\mathcal{J}^{s}(y_{T,\sigma})
		\geq & \frac{1}{2C_{obs}}\mathbb{E}(\|y_{T,\sigma}\|^2_H) -\mathbb{E}(\langle x_{1,\sigma}, y_{T,\sigma} \rangle_{H}) +\mathbb{E}(\langle x_{0,\sigma} , I^{1-\alpha}_{T-t}y_\sigma(0;y_{T,\sigma})\rangle_H)\\
		=& \frac{1}{2C_{obs}}\mathbb{E}(\|y_{T,\sigma}\|^2_H) -\mathbb{E}(\langle x_{1,\sigma}, y_{T,\sigma} \rangle_{H}) +\mathbb{E}(\langle x_{0,\sigma} , E_{\alpha}(T^{\alpha}A(\sigma))y_{T,\sigma}\rangle_H).
	\end{align*}
	Applying the inequality \eqref{Young}, we obtain
	\begin{align*}
		\mathcal{J}^{s}(y_{T,\sigma})\geq & \frac{1}{4C_{obs}}\mathbb{E}(\|y_{T,\sigma}\|^2_H) -2C_{obs}\left[\mathbb{E}(\|x_{1,\sigma}\|^2_H) -E_{\alpha}(T^{\alpha}\|A(\sigma)\|_{L^{\infty}(\mathfrak{S},\mathbb{R}^{n\times n})})\mathbb{E}(\|x_{0,\sigma}\|^2_H)\right].
	\end{align*}
	Consequently, $\mathcal{J}^s$ admits a minimum in $L^2(\mathfrak{S},H;\mu)$, and it follows from Proposition \ref{Minimizer controllability} that the system \eqref{sys} is exactly simultaneously controllable.
\end{proof}

The HUM (Hilbert Uniqueness Method) control corresponds to the control of minimal energy in the weighted space $L^{2}_{\alpha-1}(0,T;U)$.
\begin{proposition} \label{HUM Control s} Assume that system \eqref{sys} is exactly simultaneous controllable. Then, for any $x_{0,\sigma}, x_{1,\sigma}\in L^2(\mathfrak{S},H;\mu)$, the unique minimizer of the following problem:
	\[
	\min_{u} \, \|u\|_{L^2_{\alpha-1}(0,T;U)} \quad \text{subject to} \quad x_\sigma(T;x_{0,\sigma},u)=x_{1,\sigma},\quad \text{a.e},\;\sigma\in\mathfrak{S}
	\]
	is given by \eqref{control s}.
\end{proposition}
\begin{proof}
	Let $u$ be a control such that
	$x_\sigma(T;x_{0,\sigma},u)=x_{1,\sigma}$, a.e, $\sigma\in\mathfrak{S}$ and set $v=u-\hat{u}$. Assuming that 
	$v$ and $\hat{u}$ are orthogonal in $L^2_{\alpha-1}(0,T;U)$, the Pythagorean theorem yields:
	\begin{align*}
		\|u\|^2_{L^2_{\alpha-1}(0,T;U)}=\|v\|^2_{L^2_{\alpha-1}(0,T;U)} + \|\hat{u}\|^2_{L^2_{\alpha-1}(0,T;U)}.
	\end{align*}
	Consequently, $\|u\|^2_{L^2_{\alpha-1}(0,T;U)}\geq \|\hat{u}\|^2_{L^2_{\alpha-1}(0,T;U)}$ if and only if $u=\hat{u},$ a.e on $(0,T)$.\\
	Then, we conclude by showing that $\langle v,\hat{u}\rangle_{L^2_{\alpha-1}(0,T;U)}=0$. Using \eqref{control s}, we obtain
	\begin{align*}
		\langle v, \hat{u}\rangle_{L^2_{\alpha-1}(0,T;U)}&=\int_{0}^{T}(T-t)^{\alpha-1}\langle v(t),\hat{u}(t) \rangle_U dt\\
		&=\int_{0}^{T}(T-t)^{\alpha-1}\langle v(t),\mathbb{E}\left(B(\sigma)^{\top}E_{\alpha,\alpha}((T-t)^{\alpha}A(\sigma)^{\top})\hat{y}_{T,\sigma}\right) \rangle_U dt
	\end{align*}  . 
	Now, by \eqref{Duhamel formula_adj}, we can rewrite
	\begin{align*}
		\langle v, \hat{u}\rangle_{L^2_{\alpha-1}(0,T;U)}
		&=\int_{0}^{T}\langle v(t),\mathbb{E}(B(\sigma)^{\top}y_{\sigma}(t;\hat{y}_{T,\sigma})) \rangle_U dt \\
		&=\int_{0}^{T}\langle u(t),\mathbb{E}(B(\sigma)^{\top}y_{\sigma}(t;\hat{y}_{T,\sigma})) \rangle_U dt-\int_{0}^{T}\langle \hat{u}(t),\mathbb{E}(B(\sigma)^{\top}y_{\sigma}(t;\hat{y}_{T,\sigma})) \rangle_U dt.
	\end{align*}
	Finally, applying \eqref{First condition} separately to $u$ and $\hat{u}$
	yields $\langle v, \hat{u}\rangle_{L^2_{\alpha-1}(0,T;U)}=0$.
\end{proof}
\subsection{Average controllability}
\begin{Notation}
	Let $T>0$ and define, for any initial state $x_{0,\sigma}\in L^{1}(\mathfrak{S},H;\mu)$, the set of the average of the reachable states:
	\begin{align*}
		AR(T; x_{0,\sigma}):=\{\mathbb{E}(x_\sigma(T;x_{0,\sigma}; u))\colon u\in L^{\infty}(0,T;U)\},
	\end{align*}
	where, for almost every $\sigma \in \mathfrak{S}$, $x_{\sigma}(\cdot; x_{0,\sigma}; u)$ denotes the unique solution of \eqref{sys} associated with the initial data $x_{0,\sigma}$ and the control $u$.
\end{Notation}
\begin{remark}
	Proposition \ref{Well-Posedness} implies that $AR(T, x_{0,\sigma})$ is well-defined and is an affine subspace of $H$, since both the system \eqref{sys} and the expectation are linear.
\end{remark}
There are different notions of controllability in the average sense.
\begin{definition} \label{notion_controllability}
	Let $T>0$.
	\begin{itemize}
		\item System \eqref{sys} is exactly controllable in average (in time T)
		if, for every initial data $x_{0,\sigma}\in L^{2}(\mathfrak{S},H;\mu)$, the set of the average of the reachable states
		$AR(T, x_{0,\sigma})$ coincides with $H$. That is, for every $x_{0,\sigma}\in L^{1}(\mathfrak{S},H;\mu)$ and every $x_1\in H$, there exists $u\in L^{\infty}(0,T;U)$ such that:
		\begin{align}
			\mathbb{E}(x_\sigma(T;x_{0,\sigma};u))=x_{1}. \label{exactly average controllable}
		\end{align}
		\item System \eqref{sys} is approximately controllable in average (in time T)
		if, for every initial data $x_{0,\sigma}\in L^{2}(\mathfrak{S},H;\mu)$, the set of the average of the reachable states
		$AR(T, x_{0,\sigma})$ is dense in $H$. That is, for every $x_{0,\sigma}\in L^{1}(\mathfrak{S},H;\mu)$, $x_1\in H$ and every $\varepsilon>0$, there exists $u\in L^{\infty}(0,T;U)$ such that:
		\begin{align}
			\|\mathbb{E}(x_\sigma(T;x_{0,\sigma};u))-x_{1}\|_{H}<\varepsilon. \label{approximately average controllable}
		\end{align}
	\end{itemize}
\end{definition}
\begin{remark}
	Obviously, each notion of simultaneous controllability implies the corresponding notion in the average sense.
\end{remark}
\begin{remark}\label{coincide}
	In finite-dimensional spaces, affine subspaces are closed, so exact and approximate controllability in the average sense coincide.
\end{remark}

Let us first derive a necessary and sufficient condition for the exact controllability in average of \eqref{sys}.
The main idea is to characterize \eqref{exactly average controllable} through orthogonality in the space $H$:
\begin{align*}
	\langle \mathbb{E}(x_\sigma(T;x_{0,\sigma};u)), y_T\rangle_{H}=\langle x_{1}, y_{T}\rangle_{H} \quad \forall y_{T}\in H.
\end{align*}
Consequently, the duality relation \eqref{relation duality} yields the following characterization.
\begin{lemma}
	The system \eqref{sys} is driven from the initial condition $x_{0,\sigma} \in L^{1}(\mathfrak{S},H;\mu)$ to the final condition $x_1 \in H$ in time $T$ in the sense of average controllability by a control $u \in L^{\infty}(0,T;U)$, i.e., $\mathbb{E}(x_\sigma(T;x_0;u))=x_1$, if and only if
	\begin{align}
		\int_{0}^{T}\langle u(t),\mathbb{E}(B(\sigma)^{\top}y_{\sigma}(t;y_T)) \rangle_U dt=\langle x_1 , y_T \rangle_H - \mathbb{E}\left(\langle x_{0,\sigma} , I^{1-\alpha}_{T-t}y_{\sigma}(0;y_T) \rangle_H\right), \label{First condition} 
	\end{align}
	for all $y_T\in H$.
\end{lemma}
Relation \eqref{First condition} may be interpreted as an optimality condition of the critical points of the functional $\mathcal{J}^a: H\longrightarrow\mathbb{R}$:
\begin{align*}
	\mathcal{J}^a(y_T)=\frac{1}{2}\int_{0}^{T}|\mathbb{E}(B(\sigma)^{\top}y_\sigma(t;y_T))|^2_U \rho_{\alpha}(t)dt-\langle x_1 , y_T \rangle_H + \mathbb{E}\left(\langle x_{0,\sigma} , I^{1-\alpha}_{T-t}y_\sigma(0;y_T) \rangle_H\right),
\end{align*}
where $\rho_{\alpha}(t):=(T-t)^{1-\alpha},\; t\in [0,T]$.
\par 
The following two results are similar to their counterparts in the case of simultaneous controllability, and their proofs are therefore omitted.
\begin{proposition}
	Let $x_0\in L^1(\mathfrak{S},H;\mu)$, $x_1\in H$ and suppose that $\hat{y}_T\in H$ is a minimizer of $\mathcal{J}^a$. Then, for the control 
	\begin{align}
		u_a(t)=\mathbb{E}\left(B(\sigma)^{\top}E_{\alpha,\alpha}((T-t)^{\alpha}A(\sigma))\right)\hat{y}_{T} \label{control a}
	\end{align}
	we have that:
	\begin{align*}
		\mathbb{E}(x_\sigma(T;x_0;u_a))=x_1.
	\end{align*}
	
\end{proposition}

\begin{theorem}
	The following statements are equivalent:
	\begin{itemize}
		\item System \eqref{sys} fulfills the exact average controllability.
		\item System \eqref{adj_sys} fulfills the
		exact average observability:
		\begin{align}
			\exists C_{oba}>0,\; \|y_{T}\|^{2}_H\leq C_{oba} \int_{0}^{T}(T-t)^{1-\alpha}\| \mathbb{E}(B(\sigma)^{\top}y_{\sigma}(t;y_{T}))\|^{2}_{U}dt, \label{exact average observability}
		\end{align}
		for all $y_{T}\in H$.
	\end{itemize}
\end{theorem}
\begin{remark}
	The main difference in the variational approach between exact controllability and average controllability lies in the choice of test elements: in the latter, they are taken in the state space $H$, which is natural since the expectation is $\sigma$-independent.
\end{remark}
\begin{remark}
	Simultaneous controllability can also be investigated by deriving it from averaged controllability through a variational approach combined with a penalization method, as presented in Theorem 4.1 of \cite{loheac2016averaged}.
\end{remark}
\subsubsection{Average Gramian matrix}
The average controllability Gramian matrix can serve as a criterion for designing an open-loop control signal.
\begin{definition}
	The average Gramian matrix associated with system \eqref{sys} is defined by
	\begin{align}
		G_T^\alpha:=\int_{0}^{T}(T-s)^{\alpha-1}\mathbb{E}(E_{\alpha,\alpha}((T-s)^{\alpha}A(\sigma))B(\sigma))\mathbb{E}(B(\sigma)^{\top}E_{\alpha,\alpha}((T-s)^{\alpha}A(\sigma)^{\top}))ds. \label{Gramian matrix}
	\end{align}
\end{definition}
\begin{remark}
	Obviously, $G_T^\alpha\in\mathbb{R}^{n\times n}$ is symmetric and positive semidefinite matrix; that is,
	\begin{align}
		\langle G_T^\alpha y, y\rangle_H=\int_{0}^{T}(T-s)^{\alpha-1}\|\mathbb{E}(B(\sigma)^{\top}E_{\alpha,\alpha}((T-s)^{\alpha}A(\sigma)^{\top}))y\|^{2}_U ds\geq 0 \quad \forall y\in H. \label{GMP}
	\end{align}
	Using the spectral theorem, we deduce that $G_T^\alpha$ is diagonalizable in $H$ with spectrum 
	$$0\leq\lambda_1 \leq \cdots\leq \lambda_n.$$
	Furthermore, by considering an orthonormal basis of eigenvectors of $G_T^\alpha$, we obtain
	\begin{align}
		\lambda_1\|y\|^2_H \leq \langle G_T^\alpha y, y\rangle_H \leq \lambda_n \|y\|^2_H\quad \forall y\in H. \label{SPG}
	\end{align}
\end{remark}
\begin{remark}
	Based on the Duhamel formulas for the primal and adjoint systems, the Gramian matrix can be viewed through duality as follows
	\begin{align*}
		G_T^\alpha y_T=\mathbb{E}(x_{\sigma}(T;0;(T-t)^{1-\alpha}\mathbb{E}(B(\sigma)^{\top}y_{\sigma}(t;y_T))))\quad \forall y_T\in H.
	\end{align*}
\end{remark}
We now state the result used to design an open-loop control signal based on the average controllability Gramian matrix.
\begin{theorem} \label{open-loop control}
	The system \eqref{sys} is exactly average controllable if and only if the average Gramian matrix $G_T^\alpha$ is invertible. Moreover, for any $x_{0,\sigma}\in L^1(\mathfrak{S},H;\mu)$, 
	$x_1\in H$, the average associated to the following open-loop control
	\begin{align}
		\hat{u}(t)=\mathbb{E}\left(B(\sigma)^{\top}E_{\alpha,\alpha}((T-t)^{\alpha}A(\sigma)^{\top}) \right)(G_T^\alpha)^{-1}\left(x_1- \mathbb{E}(E_{\alpha}(T^\alpha A(\sigma))x_{0,\sigma})\right). \label{control formula}
	\end{align}
	satisfies $\mathbb{E}(x_\sigma(T;x_{0,\sigma};\hat{u}))=x_1$.
\end{theorem}
\begin{proof}  
	Since $G_T^\alpha$ is positive semidefinite matrix, it is invertible if and only if its smallest eigenvalue $\lambda_1$ is strictly positive. Based on \eqref{GMP}, \eqref{SPG} and \eqref{Duhamel formula_adj}, we deduce that $G_T^\alpha$ is invertible  if and only if the exact average observability \eqref{exact average observability} holds. This completes the proof of equivalence.
	Let $x_{0,\sigma}\in L^1(\mathfrak{S},H;\mu)$, 
	$x_1\in H$, and with $\hat{u}$ defined as in \eqref{control formula}, one can verify that $\mathbb{E}(x_\sigma(T;x_{0,\sigma}; \hat{u}))=x_1$. Indeed, considering this control along with the definition of the average Gramian matrix \eqref{Gramian matrix}, we obtain:
	\begin{align*}
		\int_{0}^{T}(T-s)^{\alpha-1}\mathbb{E}\left[E_{\alpha,\alpha}((T-s)^{\alpha}A(\sigma))B(\sigma)\right]\hat{u}(s)ds
		&=G_T^\alpha(G_T^\alpha)^{-1}\left(x_1- \mathbb{E}\left(E_{\alpha}(T^\alpha A(\sigma))x_{0,\sigma}\right)\right)\\
		&=x_1- \mathbb{E}\left(E_{\alpha}(T^\alpha A(\sigma)) x_{0,\sigma}\right).
	\end{align*}
	Using the Duhamel formula \eqref{Duhamel formula}, it follows that
	\begin{align*}
		\mathbb{E}(x_\sigma(T;x_{0,\sigma};\hat{u}))
		&=\mathbb{E}\left(E_{\alpha}(T^\alpha A(\sigma))x_{0,\sigma}\right) + \int_{0}^{T}(T-s)^{\alpha-1}\mathbb{E}\left[E_{\alpha,\alpha}((T-s)^{\alpha}A(\sigma))B(\sigma)\right]\hat{u}(s)ds\\
		&=x_1.
	\end{align*}
\end{proof}
\begin{remark}
	In the deterministic case $A(\sigma)=A$,  $B(\sigma)=B$ and integer-order fractional case, i.e., when 
	$\alpha=1$, we obtain:
	\begin{align*}
		G_T^1 =\int_{0}^{T}e^{(T-s)A}BB^{\top} e^{(T-s)A^\top}ds.
	\end{align*}
	We then find the usual Gramian matrix. 
\end{remark}
The following proposition characterizes the HUM control in the averaged sense, with a proof similar to that of Proposition \ref{HUM Control s}.
\begin{proposition} \label{HUM Control a}
	Assume that system \eqref{sys} is exactly average controllable. Then, for any $x_{0,\sigma}\in L^1(\mathfrak{S},H;\mu)$ and 
	$x_1\in H$, the unique minimizer of the following problem:
	\[
	\min_{u} \, \|u\|_{L^2_{\alpha-1}(0,T;U)} \quad \text{subject to} \quad \mathbb{E}(x_\sigma(T;x_{0,\sigma},u)=x_1
	\]
	is given by \eqref{control a}.
\end{proposition}
\begin{remark}
	By the uniqueness of the HUM control, the control obtained via the variational approach in \eqref{control a} coincides almost everywhere with the control derived from the averaged Gramian matrix in \eqref{control formula}.
\end{remark}
\subsubsection{Average Kalman Rank Condition}
The classical characterization of exact controllability for finite-dimensional linear systems via the Kalman rank condition is due to R. E. Kalman \cite{kalman1960contributions}. In particular, it establishes that the control time is irrelevant. The following theorem extends this result to finite-dimensional linear fractional systems, in the sense of average controllability, as considered in \cite{zuazua2014averaged} for systems with integer-order derivatives.
\begin{theorem} \label{Average Kalman Rank}
	System \eqref{sys} is exactly controllable in average if and only if
	\begin{align}
		\mbox{rank}\left[ \mathbb{E}(A^{k}B) \colon k\geq 0 \right]=n. \label{Kalman rank condition}
	\end{align}
	%	
	%	 $\mbox{rank}(\mathfrak{K})=n$, where
	%	\begin{align}
		%		\mathfrak{K}=:\left[ \mathbb{E}(B) \;\; \mathbb{E}(AB)\; \cdots\; \mathbb{E}(A^{n-1}B)\right]. \label{Kalman matrix}
		%	\end{align}
\end{theorem}

\proof
The formula \eqref{GMP} can be rewritten as follows:
\begin{align*}
	\langle G_T^\alpha y, y\rangle_H
	&=\int_{0}^{T}s^{\alpha-1}\|f_{y}(s)\|^{2}_U ds\quad \forall y\in H,
\end{align*}
where
\begin{align}
	f_y(s)=\mathbb{E}(B^{\top}E_{\alpha,\alpha}(s^{\alpha}A^{\top}))y=\sum_{k=0}^{\infty}\frac{s^{\alpha k}}{\Gamma(\alpha k+\alpha)}\mathbb{E}\left(B^{\top}(A^{\top})^{k}\right)y. \label{FY}
\end{align}
Using the orthogonality relation \eqref{orthogonality relation}, we obtain
\begin{align*}
	\mbox{rank}\left[ \mathbb{E}(A^{k}B) \colon k\geq 0 \right]<n 
%	&\Longleftrightarrow \mbox{Ker}(\mathfrak{K}^{\top})\neq \{0\} \\
	& \Longleftrightarrow \exists y_0\in H\setminus\{0\}\colon \mathbb{E}(B^{\top}(A^{\top})^{k})y_{0}=0\quad \forall k\geq 0.
\end{align*}
%	By the Cayley-Hamilton theorem, for any $k\in \mathbb{N}$, the matrix $(A^{\top})^{k}$ can be written as a linear combination of $I, A^{\top},\cdots, (A^{\top})^{n-1}$. Consequently, we obtain
%	\begin{align*}
	%		\mbox{rank}(\mathfrak{K})<n 
	%		& \Longleftrightarrow \exists y_0\in H\setminus\{0\}\colon \mathbb{E}(B^{\top}(A^{\top})^{k})y_{0}=0\quad \forall k\in\mathbb{N}.
	%	\end{align*}
By the uniqueness of the power series expansion, it follows that
\begin{align*}
	\mbox{rank}\left[ \mathbb{E}(A^{k}B) \colon k\geq 0 \right]<n 
	& \Longleftrightarrow \exists y_0\in H\setminus\{0\}\colon \sum_{k=0}^{\infty}\frac{t^{k}}{\Gamma(\alpha k+\alpha)}\mathbb{E}\left(B^{\top}(A^{\top})^{k}\right)y_0=0\quad \forall t\in (0,T^{\alpha}).
\end{align*}
Since \(\alpha>0\), the mapping \(s\mapsto t=s^\alpha\) is a bijection from \((0,T)\) onto \((0,T^\alpha)\). Hence by \eqref{FY}, we obtain
\begin{align*}
	\mbox{rank}\left[ \mathbb{E}(A^{k}B) \colon k\geq 0 \right]<n 
	& \Longleftrightarrow \exists y_0\in H\setminus\{0\}\colon f_{y_0}(s)=0\quad \forall s\in (0,T)\\
	& \Longleftrightarrow \exists y_0\in H\setminus\{0\}\colon \langle G_T^\alpha y_0, y_0 \rangle_{H}=\int_{0}^{T}s^{\alpha-1}\|f_{y_0}(s)\|^{2}_U ds=0.
\end{align*}
Since the matrix $G_T^\alpha$ is symmetric and positive semidefinite, then 
\begin{align*}
	\mbox{rank}\left[ \mathbb{E}(A^{k}B) \colon k\geq 0 \right]<n 
	&  \Longleftrightarrow G_T^\alpha\; \mbox{is not invertible}\\
	&  \Longleftrightarrow \mbox{The system}\;\eqref{sys}\;\mbox{is not controllable in average}.
\end{align*}

The averaged rank condition can be further simplified in the particular case where all the matrices \(A(\sigma)\) and \(B(\sigma)\) are scalar multiples of fixed matrices \(A\) and \(B\), respectively; that is,

$$
A(\sigma)=p(\sigma)A,\qquad B(\sigma)=q(\sigma)B.
$$
\begin{theorem} \label{Average Kalman Rank simple}
	Assume that $A(\sigma)=p(\sigma)A$ and $B(\sigma)=q(\sigma)B$ with $A$ and $B$ are constant matrices, and $p$ and $q$ are real random variables on $(\mathfrak{S},\mathcal{F},\mu)$ such that $\mathbb{E}(pq^{k})\neq 0,$ for all $k=0,1,\cdots, n-1$. Then, system \eqref{sys} is exactly controllable in average if and only if 
	\begin{align*}
		\mbox{rank}\left[ B \;\; AB\; \cdots\; A^{n-1}B\right]=n.
	\end{align*}
\end{theorem}
\begin{proof}
	Using the orthogonality relation \eqref{orthogonality relation}, and the Cayley-Hamilton theorem, for any $k\in \mathbb{N}$, the matrix $(A^{\top})^{k}$ can be written as a linear combination of $I, A^{\top},\cdots, (A^{\top})^{n-1}$. Consequently, we obtain
	\begin{align*}
	\mbox{rank}\left[ B \;\; AB\; \cdots\; A^{n-1}B\right]<n 
		& \Longleftrightarrow \exists y_0\in H\setminus\{0\}\colon B^{\top}(A^{\top})^{k}y_{0}=0\quad \forall k\in\mathbb{N}.
	\end{align*}
	By the uniqueness of the power series expansion and $\mathbb{E}(pq^{k})\neq 0,$ for all $k=0,1,\cdots, n-1$, it follows that
	\begin{align*}
		\mbox{rank}\left[ B \;\; AB\; \cdots\; A^{n-1}B\right]<n 
		& \Longleftrightarrow \exists y_0\in H\setminus\{0\}\colon \sum_{k=0}^{\infty}\frac{t^{k}}{\Gamma(\alpha k+\alpha)}\mathbb{E}\left(pq^k\right)B^{\top}(A^{\top})^{k}y_0=0\quad \forall t\in (0,T^{\alpha}).
	\end{align*}
	Since \(\alpha>0\), the mapping \(s\mapsto t=s^\alpha\) is a bijection from \((0,T)\) onto \((0,T^\alpha)\). Hence by \eqref{FY}, we obtain
	\begin{align*}
		\mbox{rank}\left[ B \;\; AB\; \cdots\; A^{n-1}B\right]<n 
		& \Longleftrightarrow \exists y_0\in H\setminus\{0\}\colon f_{y_0}(s)=0\quad \forall s\in (0,T)\\
		& \Longleftrightarrow \exists y_0\in H\setminus\{0\}\colon \langle G_T^\alpha y_0, y_0 \rangle_{H}=\int_{0}^{T}s^{\alpha-1}\|f_{y_0}(s)\|^{2}_U ds=0,
	\end{align*}
	where $f_y$ is defined in \eqref{FY}.
	Since the matrix $G_T^\alpha$ is symmetric and positive semidefinite, then 
	\begin{align*}
		\mbox{rank}\left[ B \;\; AB\; \cdots\; A^{n-1}B\right]<n 
		&  \Longleftrightarrow G_T^\alpha\; \mbox{is not invertible}\\
		&  \Longleftrightarrow \mbox{The system}\;\eqref{sys}\;\mbox{is not controllable in average}.
	\end{align*}
\end{proof}
\begin{example}
	For 
	$A(\sigma)=p(\sigma)\begin{pmatrix}
		0 & -1\\ 1 & 0
	\end{pmatrix}$ and $B=\begin{pmatrix}
		1 \\ 0
	\end{pmatrix}$, where $p$ a real random variable on $(\mathfrak{S},\mathcal{F},\mu)$ such that $\mathbb{E}(p)\neq 0$. Then, the average controllability holds.
\end{example}

\section{Some numerical results and experiments} \label{Sec5} 
\subsection{Algorithm for Computing the HUM Control}

Regarding the numerical aspects of HUM control, we refer the reader to \cite{boyer2013penalised, glowinskiexact} for detailed developments in the context of partial differential equations. In \cite{boutaayamou2025null}, this approach is applied to a heat equation coupled with an ordinary differential equation. For numerical investigations related to average and ensemble controllability in finite-dimensional systems, we refer to \cite{lazar2022control}.
\par 
Let $x_{0,\sigma} \in L^1(\mathfrak{S}, H; \mu)$ be a given initial state to be controlled, and let $x_1 \in H$ be a desired target state. In order to numerically approximate the open-loop HUM control $\hat{u}$ such that the associated average of the states at time $T$ approaches $x_1$, we make use of Theorem \ref{open-loop control}, which can be alternatively stated as follows: 
\begin{eqnarray}
	\hat{u}(t)=\mathbb{E}\left(B(\sigma)^{\top}E_{\alpha,\alpha}((T-t)^{\alpha}A(\sigma)^{\top})\right)\hat{y}_{T},  \label{control}
\end{eqnarray}
where $\hat{y}_{T}$ solves the following linear system:
\begin{eqnarray}
	G_T^\alpha\hat{y}_{T}=x_1-\mathbb{E}(x_\sigma(T,x_{0,\sigma};0)). \label{LS}
\end{eqnarray}
Since $G_T^\alpha$ is symmetric and positive definite, to resolve the linear system \eqref{LS} we proceed numerically using the Conjugate Gradient (CG) method in Algorithm \ref{algo1}, which is an efficient iterative algorithm well-suited for such matrices. CG only requires matrix-vector products and avoids the explicit computation of $(G_T^\alpha)^{-1}$, which significantly reduces computational cost and improves numerical stability. In particular, we avoid computing the inverse of the Gramian matrix, which would be both expensive and potentially ill-conditioned, especially in high dimensions.
\begin{algorithm}[H]
	\caption{Conjugate Gradient Method (CG)} \label{algo1}
	\begin{algorithmic}[1]
		\State \textbf{Input:} An initial state to be controlled \( x_{0,\sigma}\in H \), a desired target state $x_1$, avergae Gramian matrix \( G_T^\alpha \), initial guess \( y_0 \), tolerance \( \text{tol} \), maximum iterations \(k_{\text{max}}\) 
		\State \textbf{Initialize:} 
		\State $r_0=x_1-\mathbb{E}(x_\sigma(T,x_{0,\sigma};0))-G_T^\alpha y_{0}$ (an initial residual)
		\State $d_0 = r_0$ (initial descent direction)\State $k=0$
		
		\While{$||r_k|| > \text{tol}$ and $k < k_{\text{max}}$}
		\State $a_k = \frac{r_k^\top r_k}{d_k^\top G_T^\alpha d_k}$ 
		\State $y_{k+1} = y_k + a_k d_k$
		\State $r_{k+1} = r_k - a_k G_T^\alpha d_k$
		
		\If{$||r_{k+1}|| \leq \text{tol}$}
		\State \textbf{break}
		\EndIf
		
		\State $b_k = \frac{r_{k+1}^\top r_{k+1}}{r_k^\top r_k}$
		\State $d_{k+1} = r_{k+1} + b_k d_k$
		\State $k = k + 1$
		\EndWhile
		
		\State \textbf{Output:} Approximate solution $y_k$ of $\hat{y}_{T}$.
	\end{algorithmic}
\end{algorithm}
Concerning now the convergence rate in Algorithm~\ref{algo1}, we know from \cite{glowinskiexact, saad2003iterative} that
\[\|y_k -\hat{y}_T\|_{\Lambda}\leq 4\|y_0 -\hat{y}_T\|_{\Lambda}e^{-2\mathcal{C}_{GC}k},\]
where $\rho$ denotes the condition number of the bilinear functional associated with $\Lambda$, defined by $\rho=\|\Lambda\|_{_{\mathbb{R}^{n\times n}}}\|\Lambda^{-1}\|_{_{\mathbb{R}^{n\times n}}}$ and the positive constant $\mathcal{C}_{GC}$ is given by
\[\mathcal{C}_{GC}=\ln\left(\frac{\sqrt{\rho}+1}{\sqrt{\rho}-1}\right)\]
and the norm $\|\cdot\|_{\Lambda}$ is defined by $\|z\|_{\Lambda}:=\langle \Lambda(z), z\rangle_H$.

However, as seen directly from Algorithm~\ref{algo1}, this process requires computing the averaged state and adjoint equations at each iteration. To achieve this, we first generate a sample ${\sigma_1, \cdots, \sigma_M}$ of size $M$ drawn from the random parameter $\sigma$. Then, to compute the averaged state, we solve equation \eqref{sys} for each realization $\sigma_k$, $k = 1, \cdots, M$, and approximate the expectation using the classical Monte Carlo estimator $\mathbb{E}_M$:
\begin{eqnarray*}
	\mathbb{E}_{M}(x_\sigma(t; x_{0,\sigma};u))=:\frac{1}{M}\sum_{k=1}^{M}x_{\sigma_k}(t; x_{0,\sigma_k};u) \approx \mathbb{E}(x_\sigma(t; x_{0,\sigma};u))\quad \text{as} \quad M\longrightarrow \infty.
\end{eqnarray*}
More specifically, for all $t\in [0,T]$ the statistical error in the $L^2$-setting is given by (see \cite{ali2017multilevel}):
\begin{eqnarray*}
	\|\mathbb{E}(x_\sigma(t; x_{0,\sigma};u))-\mathbb{E}_{M}(x_\sigma(t; x_{0,\sigma};u))\|_{L^{2}(\mathfrak{S};H;\mu)} \leq \frac{\|x_\sigma(t; x_{0,\sigma};u)\|_{L^{2}(\mathfrak{S};H;\mu)}}{\sqrt{M}},\;\; \forall M\in \mathbb{N}\setminus \{0\}.
\end{eqnarray*}
That is, this randomized approach converges in expectation with rate $\frac{1}{2}$.
\subsection{Numerical Illustrations for the Fractional Rössler System}
The dynamics of the fractional-order Rössler system \cite{zheng2016predictive, zheng2010controlling} are governed by the following equations:
\begin{equation}\label{nonlinear Rössler}
	\begin{cases}
		\partial^{\alpha}_{t}x(t) =-y(t)-z(t),\\
		\partial^{\alpha}_{t}y(t) =x(t) + ay(t),\\
		\partial^{\alpha}_{t}z(t) =bx(t) -(c-x(t))z(t).
	\end{cases}
\end{equation}
For the parameter values $\alpha=0.97$ and $a=0.34$, $b=0.4$, $c=4.5$ the system exhibits chaotic behavior. The parameter (linear damping parameter)
$a$ controls the damping in the $y$-component of the Rössler system. 
For simplicity, we consider $a=a(\sigma)$ as a random variable, while the other parameters are fixed and study the averaged controllability of the linearized system of \eqref{nonlinear Rössler} around the equilibrium point $(0,0,0)$ with a single control input:
\begin{equation}\label{linear Rössler}
	\begin{cases}
		\partial^{\alpha}_{t}x_{\sigma}(t) =-y_{\sigma}(t)-z_{\sigma}(t),\\
		\partial^{\alpha}_{t}y(t) =x_{\sigma}(t) + a(\sigma) y_{\sigma}(t),\\
		\partial^{\alpha}_{t}z(t) =b x_{\sigma}(t) -c z_{\sigma}(t)+u(t).
	\end{cases}
\end{equation}
In this case, $A(\sigma)=\begin{pmatrix}
	0 & -1 & -1 \\ 1 & a(\sigma) & 0 \\ b & 0 & -c
\end{pmatrix}$ and $B=\begin{pmatrix}
	0 \\ 0 \\ 1
\end{pmatrix}$. A straightforward computation shows that
\[
\det\bigl(B\;\; \mathbb{E}(A)B\;\; \mathbb{E}(A^2)B\bigr)=1.
\]
Hence, the Kalman rank condition \eqref{Kalman rank condition} is satisfied.
Consequently, system \eqref{linear Rössler} is exactly average controllable. 

We now present numerical simulations for the control time \(T=2\) and the initial condition \((x_0,y_0,z_0)=(1,1,1)\). We consider two probability distributions for the uncertain parameter \(a\): a truncated Gaussian distribution and a uniform distribution.

\subsubsection{Truncated Gaussian distribution}

\begin{figure}[H]
	\centering
	\includegraphics[width=0.4\linewidth]{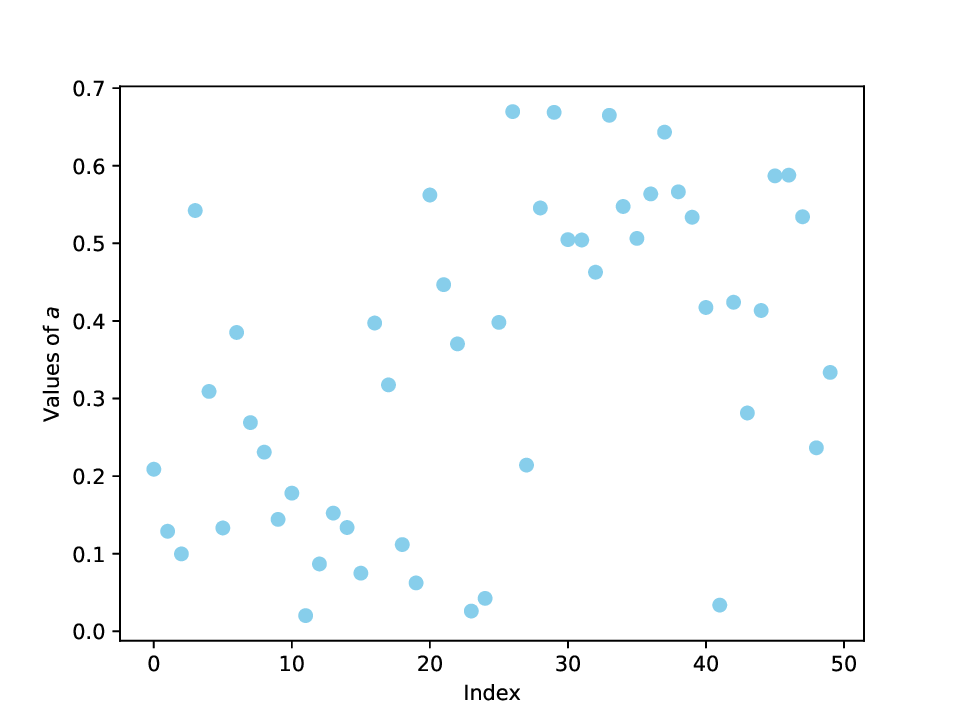}
	\caption{Sample of the truncated Gaussian distribution of the parameter \(a\).}
	\label{fig: fig1}
\end{figure}

We first consider the case in which the parameter \(a\) follows a Gaussian distribution truncated to the interval \([0,0.68]\), with mean \(0.34\) and variance \(0.5\). Figure \ref{fig: fig1} displays a sample of the resulting distribution.

\begin{figure}[H]
	\centering
	\includegraphics[width=1.1\linewidth]{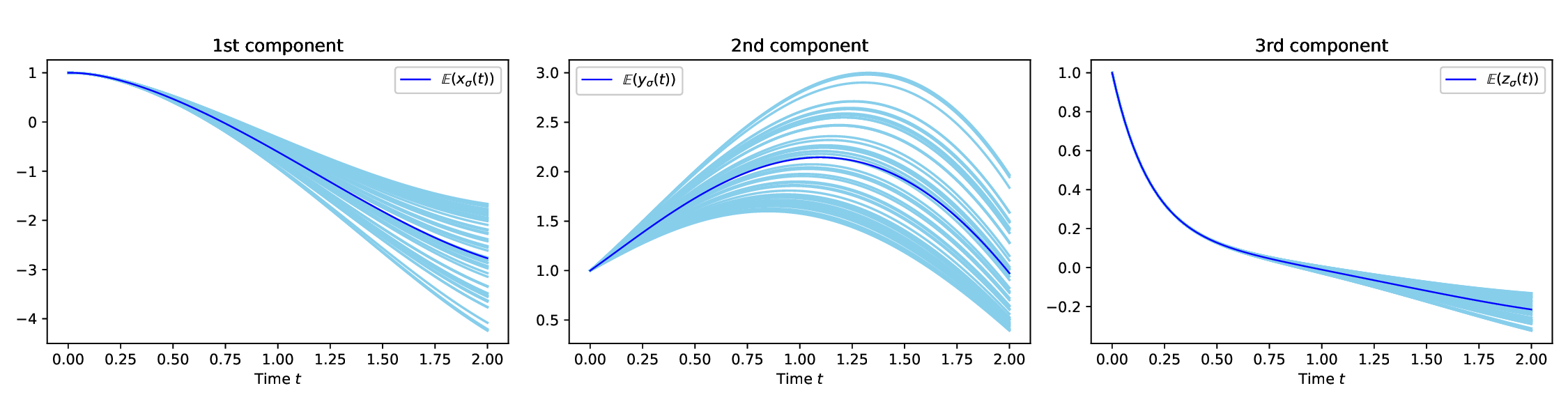}
	\caption{Uncontrolled trajectories and their ensemble average for the truncated Gaussian distribution.}
	\label{fig: fig2}
\end{figure}

Figure \ref{fig: fig2} shows the uncontrolled trajectories of system \eqref{linear Rössler} corresponding to the sampled values of \(a\), together with their ensemble average. The control is set to \(u=0\).

\begin{figure}[H]
	\centering
	\includegraphics[width=0.4\linewidth]{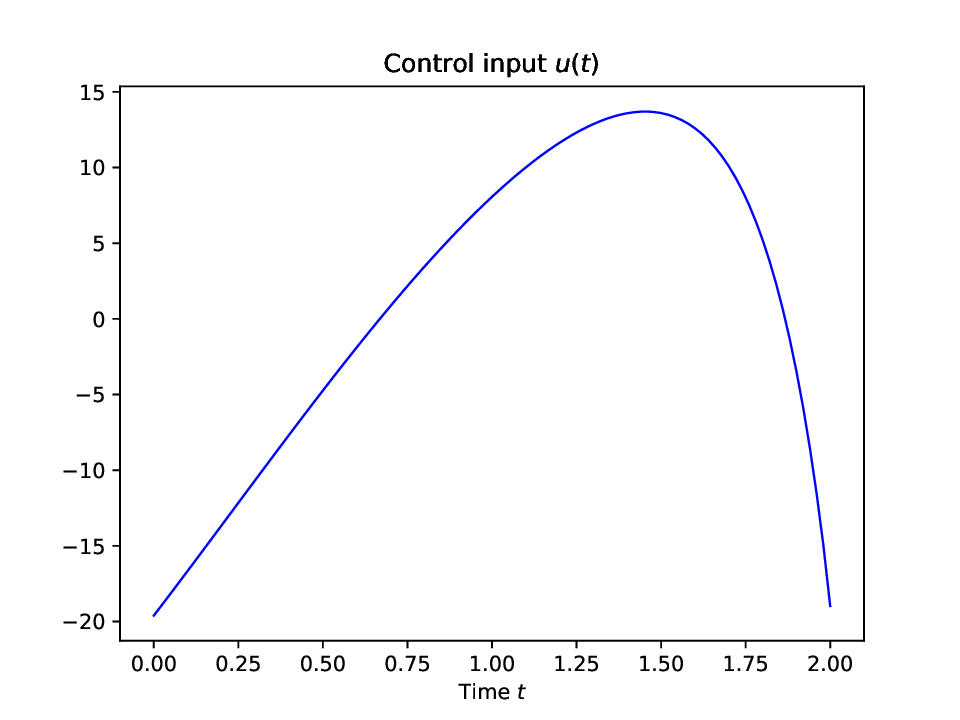}
	\caption{Numerical approximation of the HUM control for the truncated Gaussian distribution.}
	\label{fig: fig3}
\end{figure}

We next compute a numerical approximation of the Hilbert Uniqueness Method (HUM) control designed to steer the ensemble average of the system to the equilibrium point \((0,0,0)\). The resulting control function is shown in Figure \ref{fig: fig3}.

\begin{figure}[H]
	\centering
	\includegraphics[width=1.1\linewidth]{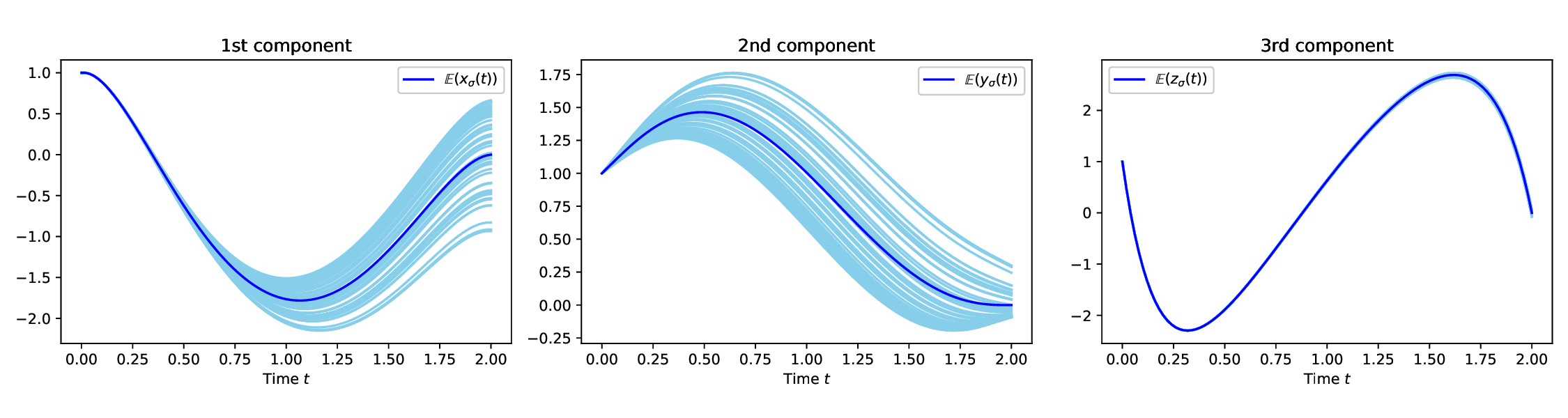}
	\caption{Controlled trajectories and their ensemble average for the truncated Gaussian distribution.}
	\label{fig: fig4}
\end{figure}

Figure \ref{fig: fig4} presents the corresponding controlled trajectories for the sampled parameter values, together with their ensemble average. The numerical algorithm converges after \(16\) iterations. At the final time \(T=2\), the computed ensemble-average state satisfies
\[
\begin{aligned}
	\mathbb{E}[x_\sigma(T)]&=-1.37\times 10^{-7},\\
	\mathbb{E}[y_\sigma(T)]&=-4.14\times 10^{-7},\\
	\mathbb{E}[z_\sigma(T)]&= 3.04\times 10^{-7}.
\end{aligned}
\]
These values indicate that the ensemble average is driven numerically very close to the target equilibrium.

\subsubsection{Uniform distribution}

\begin{figure}[H]
	\centering
	\includegraphics[width=0.4\linewidth]{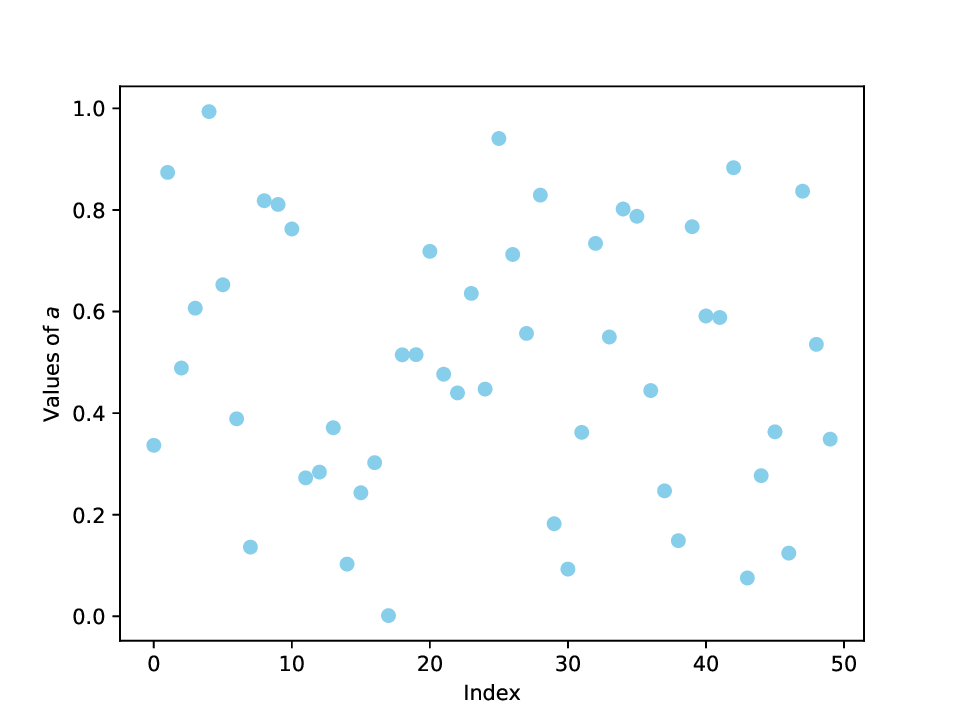}
	\caption{Sample of the uniform distribution of the parameter \(a\).}
	\label{fig: fig5}
\end{figure}

We now consider the case in which the parameter \(a\) follows a uniform distribution in $[0,1]$. Figure \ref{fig: fig5} displays a sample of this distribution.

\begin{figure}[H]
	\centering
	\includegraphics[width=1.1\linewidth]{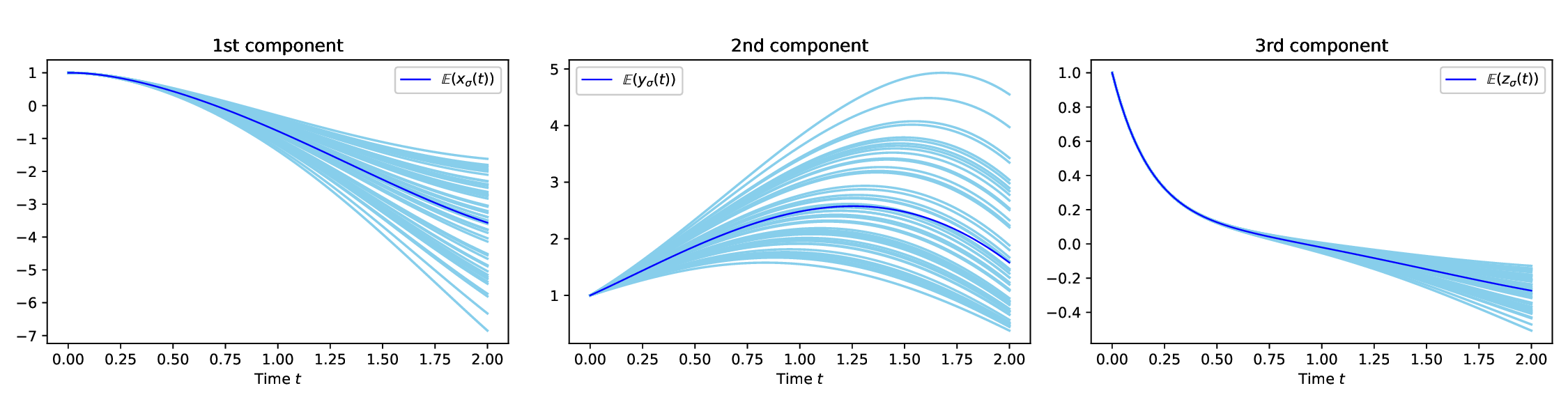}
	\caption{Uncontrolled trajectories and their ensemble average for the uniform distribution.}
	\label{fig: fig6}
\end{figure}

Figure \ref{fig: fig6} depicts the uncontrolled trajectories of system \eqref{linear Rössler} associated with the sampled values of \(a\), together with their ensemble average. As in the previous case, we take \(T=2\), \((x_0,y_0,z_0)=(1,1,1)\), and \(u=0\).

\begin{figure}[H]
	\centering
	\includegraphics[width=0.4\linewidth]{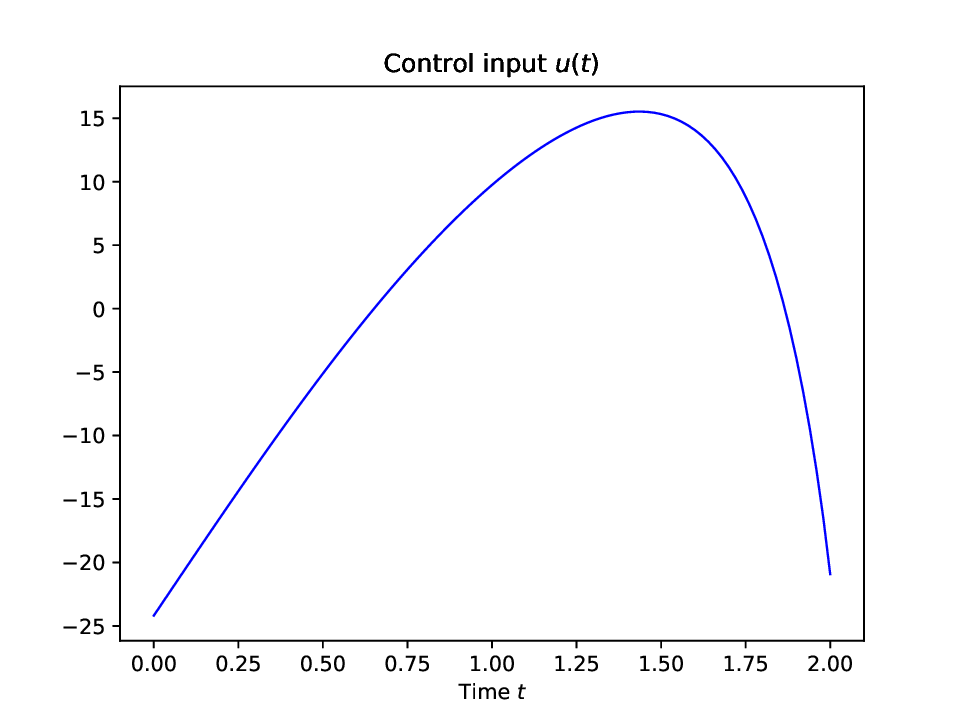}
	\caption{Numerical approximation of the HUM control for the uniform distribution.}
	\label{fig: fig7}
\end{figure}

We then compute the numerical approximation of the HUM control intended to steer the ensemble average toward the equilibrium point \((0,0,0)\). The computed control is displayed in Figure \ref{fig: fig7}.

\begin{figure}[H]
	\centering
	\includegraphics[width=1.1\linewidth]{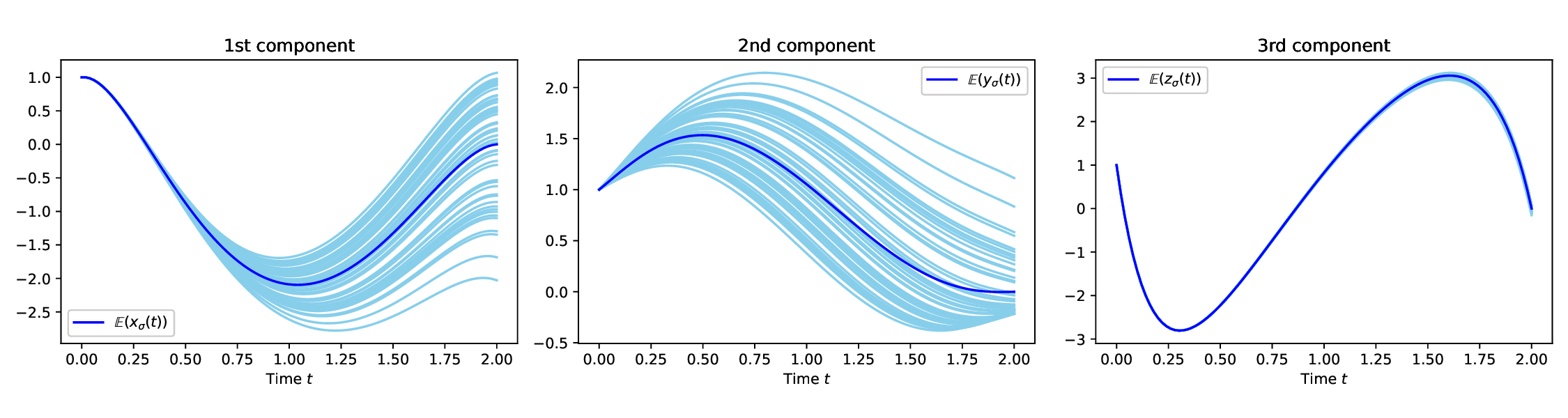}
	\caption{Controlled trajectories and their ensemble average for the uniform distribution.}
	\label{fig: fig8}
\end{figure}

Finally, Figure \ref{fig: fig8} shows the controlled trajectories corresponding to the sampled values of \(a\), along with their ensemble average. The numerical algorithm converges after \(12\) iterations. At the final time \(T=2\), the computed ensemble-average state satisfies
\[
\begin{aligned}
	\mathbb{E}[x_\sigma(T)]&=-9.82\times 10^{-8},\\
	\mathbb{E}[y_\sigma(T)]&=-3.29\times 10^{-7},\\
	\mathbb{E}[z_\sigma(T)]&=-1.16\times 10^{-7}.
\end{aligned}
\]
These residual values provide numerical evidence that the ensemble average reaches the target equilibrium with high accuracy.

\section{Conclusion}\label{Sec6}
We have investigated the average and simultaneous controllability of finite-dimensional fractional systems under uncertain parameters, using open-loop controls that are independent of the random parameter. Our results characterize average controllability through matrix conditions, specifically the Kalman rank condition and the Gramian matrix. The theoretical findings are illustrated through numerical simulations on the Rössler system.

\section*{CRediT authorship contribution statement}
\textbf{Idriss Boutaayamou}: Writing-review \& editing, Writing – original draft, Visualization, Validation, Software, Methodology, Investigation, Formal analysis, Conceptualization; \textbf{Fouad Et-Tahri}: Writing-review \& editing, Writing-original draft, Visualization, Validation, Software, Methodology, Investigation, Formal analysis, Conceptualization; \textbf{Lahcen Maniar}: Writing-review \& editing, Writing-original draft, Visualization, Validation, Software, Methodology, Investigation, Formal analysis, Conceptualization.

\section*{Data availability}
No data was used for the research described in the article.

\section*{Declaration of competing interest}
The authors declare that they have no known competing interests that could have appeared to influence the work reported in this paper.

\end{document}